\documentclass{amsart}
\usepackage{amssymb}

\def\dessous#1\sous#2{\mathrel{\mathop{\kern0pt#2}\limits_{#1}}}

\newcommand{\Z}{\mathbb Z}
\newcommand{\R}{\mathbb R}

\newcommand{\F}{\mathbb F}

\newcommand{\cc}{\mathcal{C}}
\newcommand{\A}{\mathcal{A}}
\newcommand{\X}{\mathcal{X}}
\newcommand{\I}{\mathcal{I}}
\newcommand{\el}{\mathcal{L}}
\newcommand{\G}{\mathcal{G}}
\newcommand{\B}{\mathcal{B}}
\newcommand{\f}{\mathcal{F}}
\newcommand{\E}{\mathcal{E}}
\newcommand{\h}{\mathcal{H}}
\newcommand{\M}{\mathcal{M}}
\newcommand{\n}{\mathcal{N}}
\newcommand{\ho}{\hbox{\rm Hom}}

\newcommand{\0}{/ \! \! \! 0}

\newtheorem{pro}{Proposition}[section]
\newtheorem{cor}{Corollary}[section]

\newtheorem{lem}{Lemma}[section]
\newtheorem{Th}{Theorem}[section]
\newtheorem{df}{Definition}[section]

\begin{document}

\title{The cohomology of Koszul-Vinberg algebras}

\author{Michel  NGUIFFO BOYOM}

\date{ }

\begin{abstract}
This work is devoted to an intrinsic cohomology theory of Koszul-Vinberg algebras
and their modules. Our results may be regarded as improvements of the attempt by
Albert Nijenhuis in [NA]. The relationships between the cohomology theory developed
here and some classical problems are pointed out, e.g. extensions of algebras and
modules, and deformation theory. The real Koszul-Vinberg cohomology of locally flat
manifolds is initiated. Thus regarding the idea raised by M. Gerstenhaber we can
state : The category of KV-algebras has its proper cohomology theory.
\end{abstract}
\address{D\'epartement de Math\'ematiques, Universit\'e Montpellier II}

\maketitle

\underline{\LARGE{\bf{Introduction}}}

\subsection{ } According to M. Gertenhaber, ``every restricted deformation theory generates its
proper cohomology theory", [GH]. The deformation theory of associative algebras
(resp. Lie algebras) and their modules involves the Hochschild cohomology theory of
those algebras and modules. Today there doesn't exist any standard way to associate
a proper cohomology theory to every given category of algebras and modules. The
first attempt to define a cohomology theory of Koszul-Vinberg algebras goes back to
Nijenhuis' paper [NA]. The role played by Koszul-Vinberg algebra in differential
geometry and in algebraic-analytic geometry is quite important. Furthermore, the
deformation theory of those algebras is related to the theory of quantization by
deformation; this aspect will be the purpose of a forthcoming paper $[NB5]$.

Given a Koszul-Vinberg algebra $\A$ and an $\A$-bimodule $W$, we denote by $\A_L$
the Lie algebra whose bracket operation is given by $[a, b] = ab - ba$. Then $W$
become a left $\A_L$-module, as well as the space $\el (\A, W)$ of linear maps from
$\A$ to $W$. Let $\cc(\A_L, \el(\A, W))$ be the Hochschild complex of $\el (\A,
W)$-valued cochains of $\A_L$. Nijenhuis [NA] set $H^q(\A, W) = H^{q-1} (\A_L, \el
(\A, W))$ and called it the $q^{th}$ cohomology space of the Koszul-Vinberg algebra
$\A$ with coefficients in $W$. However the above definition of A. Nijenhuis
collapses at the level $\leq 1$, (A. Nijenhuis already raised this problem).

The main purpose of the present work is to initiate an intrinsic cohomology theory
of Koszul-Vinberg algebras and their modules. We relate this cohomology theory to
classical problems such as:

$1^{\circ} \ H^0(\A, -) \leftrightarrow \A$-equivariant objects,

$2^{\circ} \ H^1(\A, -) \leftrightarrow$ Extensions of $\A$-modules,

$3^{\circ} \ H^2(\A, -) \leftrightarrow$ Extension classes of Algebras,

$4^{\circ} \ H^2(\A, -) \leftrightarrow$ Deformation theory of algebraic structures,

$5^{\circ} \ H^3(\A, -) \leftrightarrow$ Formal deformations of algebraic
structures.

\subsection{Content of the paper}

The present paper consists of three parts.

Part I is the theoretic part of the subject. Through Section 3 main definitions and
examples are given. Section 4 is devoted to the intrinsic cohomology theory of
Koszul-Vinberg algebras and their modules. The cohomology spaces $H_0$ and $H^1$ of
some important examples are computed. Section 5 and Section 6 are devoted to some
applications.

Part II is devoted to the cohomology theory of Koszul-Vinberg algebras which have
their origin in Differential Geometry and Analytical-Algebraic Geometry, [EV],
[KJL], [KV], [MJ], [NB3].

To every locally flat manifold $(M, D)$ is attached a $\Z_2$-graded Koszul-Vinberg
algebra (supporting a Lie super algebra). Those data admit some natural filtrations
leading to spectral sequences. Those spectral sequences are not studied in the
present paper.

Real Koszul-Vinberg cohomology of locally flat manifolds is defined in section 7.

In section 8 the real Koszul-Vinberg cohomology is used to examine the rigidity of
hyperbolic locally flat manifolds. The case of the cone $\R^+ \times \R$ is studied.
It is shown that hyperbolic affine structures in $\R^+ \times \R$ admit non trivial
deformations. That is a particular case of a general rigidity theorem by J. L.
Koszul [JLK].

Section 9 is devoted to the relationship between the real Koszul-Vinberg cohomology
of locally flat manifolds and their completeness. In particular every unimodular
locally flat manifolds (M, D) supports a real Koszul-Vinberg cohomology class
$[\omega] \in H^n (M, \R)$ which is represented by a D-parallel volume form. If $M$
is simply connected $[\omega]$ is an obstruction to the completeness of $(M, D)$.
Therefore if $(M, D)$ is the universal covering of a compact unimodular locally flat
manifold $(M_0, D_0)$ the real Koszul-Vinberg cohomology class $[\omega]$ may be
regarded as an alter ego of the conjecture of Markus, [FGH], [CY].

Part III also is a theoretic part. It contains an introduction to so-called left
invariant locally flat structures on groups of diffeomorphisms. The set of left
invariant locally flat structures on a finite dimensional Lie group G is
parametrized by a subset of the set of Koszul-Vinberg multiplications in the Lie
algebra of left invariant vector fields on G. On the other hand every locally flat
structure (M, D) on a manifold M gives rise to a Koszul-Vinberg algebra structure on
Diff(M). The example $(f \frac{d}{dx}) . (g \frac{d}{dx}) = fg \frac{d} {dx}$ ????
in the space of smooth vector fields on $\R$ shows that the inverse of the above
statement fails to hold. So the question to know what do left invariant locally flat
structures on $Diff(M)$ which are given by locally flat structures on M look like is
the main concern of Part III. We show that solutions to the question raised above
involve some special cohomology classes (Theorem 10.2).

\part{KV-Cohomology theory}

\section{KV-algebras}

Let F be a commutative field of characteristic zero. Let $\A$ be an algebra over F.
For elements $a$ and $b$ in $\A$ the multiplication map is denoted by $(a,b)
\longmapsto ab$. Given three elements $a$, $b$ and $c$ in $\A$ we will denote by
$(a,b,c)$ the associator of these elements, viz.
$$
(a,b,c) = (ab)c - a(bc). \leqno (1)
$$

\begin{df} An algebra $\A$ \ is called a KV-algebra, (KV for Koszul-Vinberg) if
$(a,b,c) = (b,a,c)$ for every triple $a,b,c$ in $\A$.
\end{df}

\underline{Examples of KV-algebras}.

(i) Every associative algebra is a Koszul-Vinberg algebra.

(ii) Let $F$ be the field $\R$ of real numbers and let $\A$ be the vector space
$C^\infty(\R, \R)$ of smooth functions. For $f$ and $g$ in $\A$ we set $fg = f
\frac{dg}{dx}$.

(iii) Let $M$ be a smooth manifold and let $D$ be a linear connection. Let us
suppose $D$ to be torsion free. If the curvatuve tensor of D is zero then the
multiplication map $ (a,b) \longrightarrow ab : = D_a b $ defines a KV-structure in
$\A = \X(M)$, the space of smooth vector fields on the manifold $M$.

Examples (ii) and (iii) are of infinite dimension. Finite dimensional KV-algebras
are related to the geometry of bounded domains $[KJ]_2, [VK]$, while Example (iii)
is related to affine geometry. We are going to define a chain complex that provides
a good framework for investigating the topology of hyperbolic affine manifolds
$(M,D)$.

\underline{Remark} Koszul-Vinberg algebras are also called ``left symmetric
algebras" because of the property $(a,b,c) = (b,a,c)$, see [NB3], [PA], [NA].

\begin{df} The Jacobi elements in a KV-algebra $\A$ are those elements $\xi \in
\A$ which satisfy the identity $(a,b,\xi) = 0$ for all elements a and b in $\A$.
\end{df}

The subspace of Jacobi elements in $\A$ is denoted by $J(\A)$; it is easy to verify
that $J(\A)$ is a subalgebra of $\A$ which is associative. Moreover the center
$C(\A)$  of $\A$ is contained in $J(\A)$. Indeed, if $c \in C(\A)$, and $a, b \in
\A$ then we have
$$
(ab)c - a(bc)  = (ba)c - b(ac) = c(ba) - b(ca) = (cb)a - (bc)a  =  0
$$

\section{Koszul-Vinberg modules}

Let $\A$ be a KV-algebra. We consider a vector space W with two bilinear maps :
$$
\left \{
\begin{array}{rclcl}
(i) \ \A \times W & \longrightarrow & W : (a,w) & \mapsto & aw, \\
(ii) \ W \times \A & \longrightarrow & W : (w,a) & \mapsto & wa.
\end{array} \right. \leqno (2)
$$
Now for elements a,b in $\A$ and $w \in W$ one sets
$$
(a,b,w) = (ab)w - a(bw),
$$
$$
(a,w,b) = (aw)b - a(wb),
$$
$$
(w,a,b) = (wa)b - w(ab).
$$
\begin{df} A vector space $W$ with bilinear maps (2) is called an {\bf{$\A$-KV-module}}
if the identities $(a,b,w) = (b,a,w)$ and $(a,w,b) = (w,a,b)$ hold. One says that
$W$ is a left KV-module (resp. a right KV-module) if the bilinear map (2) : (ii)
(resp. (2) : (i)) is the zero map.
\end{df}

Let us recall some classical facts. We condider a KV-algebra $\A$ and a fixed
$\A$-KV-module $W$. Let $\I$ be the bilateral ideal of $\A$ which is generated by
the associators (a,b,c). If $W$ is a right KV-module for $\A$ then the action of
$\I$ on $W$ is trivial, so that $W$ becomes a right module for the associative
algebra $\A /\I$.

\underline{Examples of KV-modules}

1) For every KV-algebra $\A$ let $|\A|$ be the underlying vector space; then $|\A|$
is a $\A$-KV-module.

2) Let $(M,D)$ be a locally flat manifold. The vector space of smooth vector fields
$\X(M)$ becomes a KV-algebra with the  multiplication $a.b = D_a b$. Let $W$ be the
vector space of smooth real valued functions $C^\infty(M, \R)$. By setting $ a.w =
<dw,a> , \ a \in \A, w \in W,$ one gets a structure of left KV-module in $W$.

We will extend the notion of Jacobi elements to KV-modules.

\begin{df} The Jacobi elements in a KV-module $W$ are those elements $w_0 \in W$
such that $(a,b,w_0) = 0$ for all $a$ and $b$ in $\A$.
\end{df}

\underline{Examples}

1) Consider the vector space $C^\infty([0,1], \R)$ of smooth functions defined in
the interval $[0,1]$. For $f$ and $g$ in $C^\infty([0,1], \R)$, one defines the
multiplication $ fg = f \frac{dg}{dx}$. We already saw that the above multiplication
endows $C^\infty([0,1], \R)$ with a KV-algebra structure. The space of Jacobi
elements in $C^\infty([0,1], \R)$ is the subspace of affine functions.

2) Consider a locally flat manifold $(M,D)$ and set as before $\A = (\X(M),D)$. Then
J($\A$) is the space of infinitesimal affine transformations of $(M,D)$. In other
words, an element $\xi$ belongs to J($\A$) if and only if $L_{\xi}D = 0, L_{\xi}$
being the Lie derivation in the direction of $\xi$.

3) Consider the KV-algebra in 2), viz $\A = (\X(M),D)$. Let $W = C^\infty(M, \R)$.
Then $W$ is a left KV-module for $\A$ : $ aw = \ < dw,a >.$ The subspace of Jacobi
elements in $W$ is the subset of smooth functions $w_0$ such that $D^2w_0 = 0$;
otherwise speaking $J(W)$ is the subspace of affine functions in the locally flat
manifold $(M,D)$.

It is useful to point out the following example. Let $\A = (\X(M),D)$ where $D$ is a
torsion free flat linear connection. Let $J(\A)$ be the subset of $J(\A)$ consisting
of those elements \ $\xi \in J(\A)$ which are geodesically complete; then $J_e(\A)$
is a finite dimensional subalgebra of the Lie algebra $J(\A)$. In particular if $M$
is compact then $J_e(\A)$ is finite dimensional; so that one can consider $J(\A)$ as
a Lie subalgebra of $\X(M)$ with the Lie bracket
$$
[a,b] = ab - ba.
$$
Therefore if $\G$ is the simply connected Lie group of the Lie algebra $J(\A)$ then
$\G$ carries a bi-invariant affine structure, $(\G, \nabla); \nabla$ is a
bi-invariant torsion free linear connection with $R_{\nabla} = 0$,  where
$R_{\nabla}$ is the curvature tensor of  $\nabla$. So the affine manifold $(M,D)$
admits $\G$ as an effective group of affine transformations ; (see [TS]) for the
relationship between the $\G$-geometry and the completeness of $(M,D)$. Of course,
the notions of KV-morphisms, sub-KV-modules and factor modules hold. In particular
the image under KV-morphism of a KV-module is a KV-submodule.

\section{KV-module of linear maps}

Let $V$ and $W$ be two KV-modules of a KV-algebra $\A$. Let $L(W,V)$ be the vector
space of linear maps from $W$ to $V$. We shall consider the following bilinear maps.
$$
\begin{array}{rclcl}
\A \times L(W,V) & \longrightarrow & L(W,V) : (a,f) & \mapsto & a.f, \\
L(W,V) \times \A & \longrightarrow & L(W,V) : (f,a) & \mapsto & f.a,
\end{array}
$$
where the linear maps a.f and f.a are defined as follows :
$$
\left \{
\begin{array}{l}
(a.f)(w) = a(f(w)) - f(aw) \\
(f.a)(w) = (f(w))a
\end{array} \right. \leqno (3)
$$

Equation (3) means that for $w \in W$ we have
$$
\begin{array}{l}
(a.f)(w) = a_v(f(w)) - f(a_ww),\\
(f.a)(w) = (f(w))a_v.
\end{array}
$$
One easily verifies that the identities $ (a,b,f) = (b,a,f)$ and $ (a,f,b)= (f,a,b)$
hold for all $a$ and $b$  in $\A$ and all $f \in L(W,V)$. Thus $L(W,V)$ is also a
KV-module for the KV-algebra $\A$. If $V$ is a left module then $L(W,V)$ is also a
left module. If $V$ and $W$ are right module then so is $L(W,V)$.

Similarly, the space $L_q(W,V)$ consisting of $q$-linear maps from $W$ to $V$
$$
(w_1 \dots w_q) \longmapsto f(w_1, \dots , w_q)
$$
becomes a KV-module under the actions
$$
\begin{array}{l}
(a.f)(w_1 \dots w_q) = a(f(w_1, \dots ,w_q)) - \sum_{j=1}^q f(a_1 \dots aw_j, \dots
w_q),\\
(f.a)(w_1, \dots w_q) = (f(w_1, \dots , w_q))a.
\end{array}
$$

\section{KV-Cohomology}

In the case when $\A$ is an associative algebra (resp. a Lie algebra) and $W$ is a
module of $\A$, the theory of Chevalley-Hochschild cohomology $H^*(\A,W)$ is well
developed. Moreover that theory provides powerful tools for bringing under control
many difficult problems such as deformations and extensions of algebraic structures,
deformation quantization, and so on.

Suppose now that $\A$ is an algebra which is neither associative nor a Lie algebra.
In general it is easy to define $\A$-modules, but it is a highly non-trivial problem
to define cohomologies of $\A$. In [NA], A. Nijenhuis tried to construct a
cohomology theory for Koszul-Vinberg algebras. The task turned out to be rather
difficult. So A. Nijenhuis only initiated a theory which is related to the
Chevalley-Hochschild cohomology of associated Lie algebras. Indeed, every KV-algebra
$\A$ supports a Lie algebra structure whose bracket operation is defined by $ [a,b]
= ab - ba.$ Let us denote by $\A_L$ the Lie algebra ($\A, [\ ,\ ])$ under the above
bracket. If $W$ is a KV-module for the KV-algebra $\A$, then it becomes a left
$\A_L$-module because of the identity $ (a,b,w) = (b,a,w).$

The Chevalley-Hochschild cohomology of $\A_L$ with coefficients in $W$ is well
defined. Moreover the linear space $L(\A_L, W)$ becomes a left $\A_L$-module as well
: $ (af)(b) = a(f(b)) - f([a,b])$ for $f \in L(\A,W)$; $a,b \in \A_L$. Thus the
cohomology space $H^*(\A_L, L(\A_L,W))$ is also well defined.

\begin{df} (A. Nijenhuis). Given a KV-algebra $\A$ and a $\A$-KV-module $W$, the $q^{th}$
cohomology space of $\A$ with coefficients in $W$ is defined as follows : $H^q(\A,W)
= H^{q-1}(\A_L, L(\A_L, W)).$
\end{df}

We make the following remarks.

1) The definition of $H^q(\A,W)$, is extrinsic to the theory of KV-algebras and
KV-modules. Indeed; following A. Nijenhuis [NA]; the vector space of q-cochains of
$\A$ with coefficients in W consists of those $q$-linear maps   $f: \A^q \longmapsto
W$ which are skew symmetric with respect to the first $q-1$ arguments and the
coboundry operator is the Lie-coboundry operator. So [NA] doesn't contain any
definition of a KV-coboundry operator.

2) A. Nijenhuis observed that his definition collapses in degree zero. This fact is
not toll-free. Indeed in many cohomology theories the cohomology classes of degree
zero often give very important informations.

\subsection{Intrinsic KV-cohomology theory}

We will now introduce an intrinsic cohomology theory for KV-algebras. We shall show
that the theory is coherent. Furthermore the cohomology space of degrees zero, one
and two will be interpreted as expected.

We start by fixing a KV-algebra $\A$ and a KV-module of $\A$, denoted by $W$.

Let $q$ be a positive integer. Let $C_q(\A,W)$ be the vector space of all $q$-linear
maps from $\A$ to $W$. We shall make use of the following notations.

Recall that $C_q(\A,W)$ is a $\A$-KV-module for the two actions of $\A$ on
$C_q(\A,W)$ :
$$
\begin{array}{rcl}
(af)(a_1\dots a_q) & = & a(f(a_1 \dots a_q)) - \sum_{j=1}^q f(a_1 \dots,aa_j;\dots a_q), \\
(fa)(a_1\dots a_q) & = & (f(a_1 \dots a_q)) a,
\end{array}
$$
where $a \in \A \ {\rm and} \ f \in C_q(\A, W)$. For each $\rho = 1,...,q$, we
denote by $e_{\rho}(a)$ the linear map from $C_q(\A,W)$ to $C_{q-1}(\A,W)$ which is
defined by
$$
(e_{\rho}(a)f)(a_1 \dots a_{q-1}) = f(a_1 \dots a_{\rho-1}, a, a_{\rho} \dots
a_{q-1}).
$$
We are going now to define the coboundary operator $\delta : C_q(\A,W)
\longrightarrow C_{q+1}(\A,W)$. Let $f \in C_q(\A,W)$ and $(a_1 \dots a_{q+1}) \in
\A^{q+1}$. Then the coboundary $\delta f \in C_{q+1}(\A,W)$ is given by the
following formula
$$
(\delta f)(a_1 \dots a_{q+1}) = \sum_{1\leq j\langle q+1}(-1)^j \{(ajf)(a_1 \dots
\hat{a}_j \dots a_{q+1}) + (e_q(a_j)(fa_{q+1}))(a_1 \dots \hat{a}_j \dots \hat{a}_{q+1})
\}. \leqno (4)
$$

\begin{lem} $\delta \circ \delta = 0$.
\end{lem}

\underline{Proof}. Formula (4) implies that $\delta f(a_1 \dots a_{q+1})$ is the sum
of $q$ terms.
$$
(-1)^j((a_jf)(a_1 \dots \hat{a}_j \dots a_{q+1}) + (f(a_1 \dots \hat{a}_j \dots a_q,
a_j))a_{q+1}).
$$
Fix $i, j$ and $k$ such that $1 \leq i < j < k < q+2.$ Since $\delta(\delta f)(a_1
\dots a_{q+2})$ is the sum of the $q+1$ terms
$$
(-1)^i((a_i \delta f)(a_1\dots \hat{a}_i \dots a_{q+2}) + (\delta f(a_1 \dots
\hat{a}_i \dots a_{q+1},a_i))a_{q+2}).
$$
Set
$$
\begin{array}{rcl}
\tau_{ii} & = & a_i(\delta f(a_1\dots \hat{a}_i \dots a_j \dots a_k \dots
a_{q+2})); \\
\tau_{ij} & = & \delta f(a_1\dots \hat{a}_i \dots a_ia_j \dots a_k \dots
a_{q+2}); \\
\tau_{ik} & = & \delta f(a_1\dots \hat{a}_i \dots a_j \dots a_ia_k \dots
a_{q+2}); \\
\tau_{iq+2} & = & \delta f(a_1\dots \hat{a}_i \dots a_j \dots a_k \dots a_ia_{q+2}).
\end{array}
$$
Then we have $ (a_i \delta f)(a_1 \dots \hat{a}_i \dots a_{q+2}) = \tau_{ii} -
\sum_{\rho \not=i} \tau_{i\rho}.$ To develop the $\tau_{ii}$ as well as the
$\tau_{ij}$ we shall adopt mutatis mutandis the same scheme. So let us write
$\tau_{ii} = a_i(\sum_{\rho \leq i}(-1)^\rho \Gamma_{\rho} + \sum_{i \langle \rho
\langle q+2}(-1)^{\rho-1} \Gamma_{\rho}),$ with $ \Gamma_{\rho} =
a_{\rho}(f(a_1\dots \hat{a}_i \dots \hat{a}_{\rho} \dots a_{q+2})) - \sum_s
f(a_1\dots \hat{a}_i \dots \hat{a}_{\rho} \dots a_{\rho}a_s \dots a_{q+2}) + (f(a_1
\dots \hat{a}_1 \dots \hat{a}_{\rho} \dots a_{q+1}, a_{\rho})) a_{q+2}.$ By focusing
on the indexes $i < j < k$ we may set
$$
\begin{array}{rl}
\Gamma_j & = a_j(f(a_1\dots \hat{a}_i \dots \hat{a}_j \dots a_{q+2}) + (f(a_1 \dots
\hat{a}_i \dots \hat{a}_j \dots a_{q+1},aj))a_{q+2} \\
 & - f(a_1 \dots \hat{a}_i \dots \hat{a}_j \dots a_ja_k \dots a_{q+2}) - f(a_1 \dots \hat{a}_i
\dots \hat{a}_j \dots a_ja_{q+2}) \\
 & + \ \hbox{\rm other terms.}
\end{array}
$$
Therefore we see that the quantity $\tau_{ii}$ may be written as it follows :
$$
\hbox{$\tau_{ii} = (-1)^{j-1}a_i$} \left \{
\begin{array}{l}
a_j(f(a_1 \dots \hat{a}_i \dots \hat{a}_j \dots a_{q+2})) + (f(a_1 \dots \hat{a}_i
\dots
\hat{a}_j \dots a_{q+1},a_j))a_{q+2} \\
- f(a_1 \dots \hat{a}_i \dots \hat{a}_j \dots a_ja_k \dots a_{q+2}) - f(a \dots
\hat{a}_i
\dots \hat{a}_j \dots a_ja_{q+2}) \\
+ \ \hbox{\rm other terms}
\end{array} \right \} \hbox{+ other $\Gamma_{\rho}.$} \leqno (5)
$$
Now mutatis mutandis we get
$$
\hbox{$\tau_{jj} = (-1)^ia_j$} \left \{
\begin{array}{l}
a_i(f(a_1 \dots \hat{a}_i \dots \hat{a}_j \dots a_{q+2})) + ((f(a_1 \dots \hat{a}_i
\dots
\hat{a}_j \dots a_{q+1},a_j))a_{q+2} \\
- f(a_1 \dots \hat{a}_i \dots \hat{a}_j \dots a_ia_k \dots a_{q+2}) - f(a \dots
\hat{a}_i
\dots \hat{a}_j \dots a_ia_{q+2}) \\
+ \ \hbox{\rm other terms}
\end{array} \right \} \hbox{+ other $\Gamma_{\rho}.$} \leqno (6)
$$

One has to keep in mind that in ($\delta \circ \delta f)(a_1 \dots a_{q+2})$ the
sign of the term $\tau_{ii}$ is $(-1)^i$. We are going now to develop the terms
$\tau_{ij}$. As we have done above let us focus on the fixed indexes $i < j < k$,
then we get
$$
\hbox{$\tau_{ij} = (-1)^{j-1}$} \left \{
\begin{array}{l}
a_ia_j(f(a_1 \dots \hat{a}_i \dots \hat{a}_j \dots a_{q+2})) + (f(a_1 \dots
\hat{a}_i \dots
\hat{a}_j \dots a_{q+1}, a_ia_j))a_{q+2} \\
- f(a_1 \dots \hat{a}_i \dots \hat{a}_j \dots (a_ia_j)a_k, \dots a_{q+2}) -
\hbox{\rm other terms}
\end{array} \right \} \hbox{+ other $\Gamma_{\rho};$} \leqno (7)
$$
$$
\hbox{$\tau_{ji} = (-1)^i$} \left \{
\begin{array}{l}
a_ja_i(f(a_1 \dots \hat{a}_i \dots \hat{a}_j \dots a_{q+2})) + (f(a_1 \dots
\hat{a}_i \dots
\hat{a}_j \dots a_{q+1}, a_ja_i))a_{q+2} \\
- f(a_1 \dots \hat{a}_i \dots \hat{a}_j \dots (a_ja_i)a_k, \dots a_{q+2}) -
\hbox{\rm other terms}
\end{array} \right \} \hbox{+ other $\Gamma_{\rho}.$} \leqno (8)
$$
The signs of the terms $\tau_{ij}$ and $\tau_{ji}$ in $\delta^2f(a_1 \dots a_{q+2})$
are $(-1)^{i+1}$ and $(-1)^{j+1}$ respectively.

Develop $\tau_{ik}$ and $\tau_{jk}$ by the same method as above:
$$
\hbox{$\tau_{jk} = (-1)^i$} \left \{
\begin{array}{l}
a_i(f(a_1 \dots \hat{a}_i \dots \hat{a}_j \dots a_ja_k \dots, a_{q+2})) + (f(a_1
\dots \hat{a}_i
\dots \hat{a}_j \dots a_ja_k \dots a_{q+1}, a_i))a_{q+2} \\
- f(a_1 \dots \hat{a}_i \dots \hat{a}_j \dots a_i(a_ja_k), \dots a_{q+2}) -
\hbox{\rm other terms}
\end{array} \right \} \hbox{+ other $\Gamma_{\rho}.$} \leqno (9)
$$
$$
\hbox{$\tau_{ik} = (-1)^{j-1}$} \left \{
\begin{array}{l}
a_j(f(a_1 \dots \hat{a}_i \dots \hat{a}_j \dots a_ia_k \dots, a_{q+2})) + (f(a_1
\dots \hat{a}_i
\dots \hat{a}_j \dots a_ia_k \dots a_{q+1}, a_j))a_{q+2} \\
- f(a_1 \dots \hat{a}_i \dots \hat{a}_j \dots a_j(a_ia_k), \dots a_{q+2}) -
\hbox{\rm other terms}
\end{array} \right \} \hbox{+ other $\Gamma_{\rho}.$} \leqno (10)
$$
Of course for every $l \leq i$ the expression of $\tau_{jk}$ contains the following
terms:
$$
\hbox{$(-1)^{\ell}$} \left \{
\begin{array}{l}
a_{\ell}(f(a_1 \dots \hat{a}_{\ell} \dots \hat{a}_i \dots \hat{a}_j \dots a_ja_k
\dots, a_{q+2})) + (f(a_1 \dots \hat{a}_{\ell} \dots \hat{a}_j \dots a_ja_k \dots
a_{q+1},
a_{\ell}))a_{q+2} \\
- f(a_1 \dots \hat{a}_{\ell} \dots a_{\ell}a_i \dots \hat{a}_j \dots a_ja_k \dots
a_{q+2}) - \hbox{\rm other terms}
\end{array} \right \} \leqno (11)
$$

It is useful to develop $\tau_{\rho,q+2}$ for $\rho = i \ {\rm and} \ \rho = j$ :

$$
\begin{array}{ll}
\tau_{i,q+2} & = (-1)^{j-1} \left \{
\begin{array}{l}
a_j(f(a_1 \dots \hat{a}_i \dots \hat{a}_j \dots a_ia_{q+2})) + (f(a_1 \dots
\hat{a}_i \dots
\hat{a}_j \dots a_{q+1}, a_j))a_ia_{q+2} \\
- f(a_1 \dots \hat{a}_i \dots \hat{a}_j \dots a_ja_k \dots a_{q+1}, a_ia_{q+2}) -
\hbox{\rm other terms}
\end{array} \right \}\\
 & + {\rm other} \ \Gamma_{\rho},
\end{array} \leqno (12)
$$
$$
\begin{array}{ll}
\tau_{j,q+2} & = (-1)^i \left \{
\begin{array}{l}
a_i(f(a_1 \dots \hat{a}_i \dots \hat{a}_j \dots a_ja_{q+2})) + (f(a_1 \dots
\hat{a}_i \dots
\hat{a}_j \dots a_{q+1}, a_i))a_ja_{q+2} \\
- f(a_1 \dots \hat{a}_i \dots \hat{a}_j \dots a_ia_k \dots a_ja_{q+2}) - \hbox{\rm
other terms}
\end{array} \right \}\\
& + {\rm other} \ \Gamma_{\rho}.
\end{array} \leqno (13)
$$

Apply the same ideas to express $(\delta f. a_{q+2})(\dots \hat{a}_{\rho} \dots
a_{q+1}, a_{\rho})$ :
$$
\begin{array}{ll}
(\delta f. a_{q+2})(\dots \hat{a}_i \dots a_{q+1}, a_i) & = (-1)^{j-1}
\left \{
\begin{array}{l}
(a_j(f(a_1 \dots \hat{a}_i \dots \hat{a}_j \dots a_{q+1})))a_{q+2} \\
\quad \quad + ((f(a_1 \dots \hat{a}_i \dots \hat{a}_j \dots a_{q+1}, a_j))a_i)a_{q+2} \\
- f(a_1 \dots \hat{a}_i \dots \hat{a}_j \dots a_ja_k \dots a_{q+1}, a_i))a_{q+2} -
\hbox{\rm other terms}
\end{array} \right \} \\
 & +{\rm  other} \ \Gamma_{\rho}.
\end{array} \leqno (14)
$$
$$
\begin{array}{ll}
(\delta f.a_{q+2}(\dots \hat{a} \dots a_{q+1},a_j) & = (-1)^i \left \{
\begin{array}{l}
(a_i(f(a_1 \dots \hat{a}_i \dots \hat{a}_j \dots a_{q+1},a_j)))a_{q+2}\\
\quad \quad + ((f(a_1 \dots \hat{a}_i \dots
\hat{a}_j \dots a_{q+1},a_i))a_j)a_{q+2} \\
-(f(a_1 \dots \hat{a}_i \dots \hat{a}_j \dots a_ia_k \dots a_{q+1},a_j))a_{q+2} -
{\rm other \ terms}
\end{array} \right \} \\
& + {\rm other} \ \Gamma_{\rho}.
\end{array} \leqno (15)
$$
Now we are well positioned to prove that $(\delta \circ \delta) f = 0$. Of course,
the expression of $(\delta \circ \delta) f(a_1 \dots a_1 \dots a_{q+1})$ is nothing
but
$$
\sum_{1 \leq \rho \leq q+1}(-1)^{\rho} \{(a_\rho \Gamma f)(\dots \hat{a}_{\rho}
\dots a_{q+2}) + (\delta f(a_1 \dots \hat{a}_{\rho} \dots a_{q+1}, a \rho))a_{q+2}
\}.
$$
Keeping in mind the signs of the $\tau_{\rho,s}$ if we focus on indexes $i < j < k$,
from (5), (6), (7) and (8) we deduce that $\delta \circ \delta f(a_1 \dots a_{q+2})$
contains the following expression as its additive term :
$$
\left \{
\begin{array}{l}
(-1)^{i+j-1}a_i(a_jf(a_1 \dots \hat{a}_i \dots \hat{a}_j \dots a_{q+2})) +
(-1)^{i+j}a_j(a_if(a_1 \dots
\hat{a}_i \dots \hat{a}_j \dots a_{q+2})) \\
+ (-1)^{i+j}(a_ia_j) f(a_1 \dots \hat{a}_i \dots \hat{a}_j \dots a_{q+2})) +
(-1)^{i+j-1}(a_ja_i) f(a_1 \dots \hat{a}_i \dots \hat{a}_j \dots a_{q+2}).
\end{array} \right. \leqno (*)_1
$$

The identity $(a,b,w) = (b,a,w)$ for all $w \in W$ and $a,b \in \A$ shows that
$(*)_1$ is zero.

From the expressions (5), (12), (14) and (15) one easily sees that $\delta \circ
\delta f(a_1 \dots a_{q+2})$ contains as an additive term the expression
$$
\hbox{$(*)_2 =$} \left \{
\begin{array}{l}
(-1)^{i+j}(a_i(f(a_1 \dots \hat{a}_i \dots \hat{a}_j \dots a_{q+1},a_j)))a_{q+2} +
(-1)^{i+j-1}a_i((f(a_1 \dots \hat{a}_i \dots \hat{a}_j \dots a_{q+1},a_j))a_{q+2}) \\
+ (-1)^{i+j}(f(a_1 \dots \hat{a}_i \dots \hat{a}_j \dots a_{q+1},a_j))a_ia_{q+2} +
(-1)^{i+j-1}((f(a_1 \dots \hat{a}_i \dots \hat{a}_j \dots a_{q+1},a_j))a_i)a_{q+2}
\end{array} \right.
$$
The expression $(*)_2$ above has the form $(a_i,w,a_{q+2}) - (w,a_i,a_{q+2})$. Hence
we get $(*)_2 = 0$.

Considering the expressions (7), (8), (9) and (10) we see that $\delta \circ \delta
f(a_1 \dots a_{q+2})$ has as an additive term the quantity
$$
\hbox{$(*)_3 =$} \left \{
\begin{array}{ll}
(-1)^{i+j-1}f(a_1 \dots \hat{a}_i \dots \hat{a}_j \dots (a_ia_j)a_k \dots a_{q+2})
+ (-1)^{i+j}f(a_1 \dots \hat{a}_i \dots \hat{a}_j \dots a_i(a_ja_k) \dots a_{q+2})\\
+ (-1)^{i+j} f(a_1 \dots \hat{a}_i \dots \hat{a}j \dots a_i(a_ja_k) \dots a_{q+2}) +
(-1)^{i+j-1} f(a_1 \dots \hat{a}_i \dots \hat{a}j \dots a_j(a_ia_k) \dots a_{q+2}).
\end{array} \right.
$$

The quantity $(*)_3$ is nothing but
$$
(-1)^{i+j-1}(f(a_1 \dots \hat{a}_i \dots \hat{a}_j \dots , (a_i,a_j,a_k), \dots
a_{q+2}) -f(a_1 \dots \hat{a}_i \dots \hat{a}_j \dots (a_j,a_i,a_k) \dots a_{q+2})).
$$
Since $(a_i,a_j,a_k) = (a_j,a_i,a_k)$ one gets $ (*)_3 = 0.$

From (5) and (9) one deduces that $\delta \circ \delta f(a_1 \dots a_{q+2})$
contains twice the term $a_if(a_1 \dots \hat{a}_i \dots \hat{a}_j \dots a_ja_k \dots
a_{q+2})$ with opposite signs, hence its contribution is reduced to zero.

From (9) and (14) one deduces that $\delta \circ \delta f(a_1 \dots a_{q+2})$
contains twice the term
$$
(*)_4 = (f(a_1 \dots \hat{a}_i \dots \hat{a}_j \dots a_ja_k \dots
a_{q+1},a_i))a_{q+2}
$$
with opposite signs, so its contribution is zero.

To end the proof of Lemma 4.1 it remains to examine the term of $\delta \circ \delta
f(a_1 \dots a_{q+2})$ which has the following form: $f(a_1 \dots \hat{a}_\ell \dots
a_{\ell}a_i \dots \hat{a}_j \dots a_ja_k \dots a_{q+2}).$ Those terms come from
(14). Indeed let $l < i < j < k$. By vertu of (4), $\delta \circ \delta f(a_1 \dots
a_{q+2})$ contains twice the additive term
$$
(*)_5 = f(a_1 \dots \hat{a}_{l} \dots a_{l}a_i \dots \hat{a}_j \dots a_ja_k \dots
a_{q+2}).
$$
The first of them appears as an additive term of $ -(-1)^j\delta f(a_1 \dots
a_{\ell} \dots a_i \dots \hat{a}_j \dots a_ja_k \dots a_{q+2})$ with the sign
$(-1)^{j+l}$. The second one appears as an additive term of $ -(-1)^{l} \delta f(a_1
\dots \hat{a}_{l} \dots a_{l}a_i \dots a_j \dots a_k \dots a_{q+2}) $ with the sign
$(-1)^{j+ l -1}$. Therefore $(\delta \circ \delta) f(a_1 \dots a_{q+2})$ doesn't
contain any non zero term of the form $(*)_5$.

We just examined all types of additive terms of $(\delta \circ \delta) f(a_1 \dots
a_{q+2})$. All of them cancel each other. Therefore $\delta \circ \delta f(a_1 \dots
a_{q+2}) = 0$ for any $f \in C_q(\A,W)$ and any $(a_1 \dots a_{q+2}) \in \A^{q+2}.$
That ends the proof of Lemma 4.1.0. $\triangle$

Let us set $ C(\A,W) = \bigoplus_{q \geq 1} C_q(\A,W)$. By Lemma 4.1 the coboundary
operator $\delta$ given by (4) endows $C(\A,W)$ with the structure of graded cochain
complex
$$
... \longrightarrow C_q(\A,W) \stackrel{\delta}{\longrightarrow} C_{q+1}(\A,W)
\stackrel {\delta}{\longrightarrow} C_{q+2}(\A,W) \longrightarrow ... .
$$
The $q^{th}$ cohomology $H^q(\A,W)$ of the above complex is well defined for $q >
1$. Hence
$$
H(\A,W) = \bigoplus_{q \geq 1} H^q(\A,W)
$$
is well defined.

To complete the picture we must define $C_0(\A,W)$ and $ \delta : C_0(\A,W)
\longrightarrow C_1(\A,W)$ such that $\delta \circ \delta \ (C_0(\A,W)) = 0$. Once
$\delta : C_0(\A,W) \longrightarrow C_1(\A,W)$ is defined we shall be able to define
$H^0(\A,W)$ and $H^1(\A,W)$.

\begin{df} We set $C_0(\A,W) = J(W) \ {\rm and} \ (\delta w)(a) = - \ aw \ + \ wa$, for all $a
\in \A$ and $w \in J(W)$.
\end{df}

It is easy to see that $\delta (\delta w) = 0$ if and only if $w \in J(W)$. So that
we get the total complex

$$
C_{\tau}(\A,W) = J(W) \oplus C(\A,W).
$$

\underline{Remark}. If $W$ is a right KV-module of $\A$ then $C_0(\A,W) = W$ because
$\delta (\delta w)$ will vanish for all $w \in W$. Now regarding the question raised
by Gerstenhaber we are in position to state the following result.

\begin{Th} The category of KV-algebras admits its proper cohomology theory.
\end{Th}

Since every associative algebra is also a KV-algebra we have the following corollary.

\begin{cor} The category of associative algebras admits another cohomology theory which is
different form its Hochschild cohomology theory.
\end{cor}

\underline{Some examples}

1) Let $(M,D)$ be a locally flat smooth manifold. We set $\A = (\X(M),D)$, the
KV-algebra of smooth vector fields. Then we have $ H^1(\A,\A) = 0.$ Indeed, let
$\psi$ be a cocycle of degree one, then for $a_1 \in \A$ and $a_2 \in \A$ one gets
$D_{a_1} \psi(a_2) + D_{\psi(a_1)}a_2 = \psi (D_{a_1}a_2).$ Therefore it is easy to
see that $\psi$ is a derivation of the Lie algebra $(\X(M),[\ ,\ ])$ where $[\ ,\ ]$
is the usual Lie bracket in $\X(M)$. By the vertu of a classical theorem by Takens
(every derivation of $(\X(M), [\ ,\ ])$ is an inner derivation), there exists a
smooth vector field $\xi$ such that $ \psi(a) = [\xi,a] = -a \xi + \xi a.$ Since
$\psi$ is a 1-cocycle of the chain complex $C_{\tau}(\A,\A)$, we conclude that $\xi
\in J(\A)$, so that $\psi$ is an exact cocycle of $C_{\tau}(\A)$.

2) Let $W = C^{\infty}(M,\R)$ be the vector space of smooth real valued functions
defined on a locally flat smooth manifold $(M,D)$. Then $W$ is a left KV-module for
$\A = (\X(M),D)$ where one sets
$$
aw = < dw , a >.
$$
Let $\psi : \A \longmapsto W$ be a cocycle of degree one. Then for $a,b \in \A$ one
has $ a \psi(b) - \psi(ab) = 0.$ Therefore $\psi$ is a linear map which is
$D$-parallel. So, if $w \in W$ one may write $w(a \psi (b))  = (wa) \psi (b) = \psi
((wa)b)$. Since every vector field may be locally written as $D_ab$, we see that
$\psi (wa) = w \psi(a)$ for every $w \in W$ and for every $a \in \X(M)$. Hence
$\psi$ is an usual differential 1-form on the manifold M. Since $D \psi = 0$, $\psi$
is a de Rham cocycle.

The subspace $J(W)$ consists of smooth function $w \in C^{\infty}(M,\R)$ such that
$(ab)w - a(bw) = 0$. Such functions are affine functions on the locally flat
manifold $(M,D)$. For each $w \in J(W)$, $\delta w$ is nothing but the exterior
differential of $w$. The space of 1-cocycles in $C_1(\A,W)$ consists of locally
linear closed 1-forms, and $\delta C_0(\A,W) = \delta J(W)$ consists of the
differentials of affine functions. Therefore
$$
H^1(\A,W) = \frac{[\hbox{\rm locally linear closed 1-forms}]}{d \ \hbox{\rm [affine
functions]}} $$ where $d : C^{\infty}(M,\R) \longrightarrow \Omega^1(M,\R)$ is the
de Rham differential operator. So there is a canonical linear map $ H^1(\A,W)
\longrightarrow H^1_{DR}(M,\R)$ which is an injective map.

3) Our third example is a combination of \ 1) and 2). Consider $\A = (\X(M),D) \
{\rm and} \ W = C^{\infty}(M,\R)$ as in Example 2). We define in $\A \oplus W \cong
\A \times W$ the following multiplication :
$$
(a,w)(a',w') = (aa',aw' + ww'),
$$
where $ww'$ is the usual product of two real valued functions. If $(a'',w'') \in
J(\A \oplus W)$ then $w'' = 0 \ {\rm and} \ a \in J(\A)$. So that $J(\A \oplus W)
\simeq J(\A).$ Let $(a'',o) \in J(\A \oplus W) = C_0(\A \oplus W, \A \oplus W)$ then
$(\delta(a'',o))(a,w) = (-aa'' + a''a\ ,a''w)$. So $\delta(a'',o) = 0$ if and only
if $a'' = 0$.

Hence the boundary map $\delta : J(\A) \longrightarrow C_1(\A \oplus W, \A \oplus
W)$ is an injective map.

Let $\theta : \A \times w \longrightarrow \A \times w$ be an 1-cocycle; for $(a,w)
\in \A \oplus w$ let us set $\theta (a,w) = (\phi(a,w), \psi(a,w))$. Then for $(a,w)
\ {\rm and} \ (a',w') \ {\rm in} \ \A \times W$ one has

$$
\begin{array}{rl}
- \delta \theta ((a,w),(a',w')) & = (a,w)(\phi(a',w') , \psi(a',w')) - (\phi(aa',aw' +ww'),
\psi(aa',aw' +ww')) \\
 & +(\phi(a,w), \psi(a,w))(a',w').
\end{array}
$$

Therefore the equation $\delta \theta = 0$ gives rise to the following system of
equations

$$
\begin{array}{l}
a \phi(a',w') - \phi(aa',aw' + ww') + \phi(a,w).a' = 0 \\
a \psi(a'w') + w \psi(a'w') - \psi(aa',aw' +ww') + \phi(a,w).w' + \psi(a,w)w' = 0 .
\end{array}
$$

Since we may put $ \theta(a,w) = \theta(a,0) + \theta(0,w)$, the identity satisfied
by $\phi$ shows that $\phi(a,w)$ doesn't depend on $w$. Thus by putting $ \phi(a,o)
= \phi(a)$, we see that $ a \phi(a') + \phi(a)a' = \phi(aa').$

By the vertu of a theorem by Takens stating that each derivation of the Lie algebra
$(\X(M), [\ ,\ ])$ has the form $\dot{a} \longmapsto [\xi_0,a]$, there exists $\xi
\in \A$ such that $ \phi(a) = [\xi,a] = -a \xi + \xi a $. Because $\theta$ is an
1-cocycle we must have $\xi \in J(\A).$

Let us examine the $W$-component of $\theta$, say $\psi(a,w)$. Set $ \psi(a,w) =
\psi(a,o) + \psi(o,w) = \lambda(a) + \mu(w).$ We know that $\phi(a,w) = [\xi,a], \xi
\in J(\A)$. So that $\delta \theta = 0$ gives us
$$
\begin{array}{l}
a(\lambda(a') + \mu(w')) + w(\lambda(a') + \mu(w')) - \lambda(aa') - \mu(aw' + ww')\\
+ [\xi,a].w' + (\lambda(a) + \mu(w))w' = 0 .
\end{array}
$$

Thus one sees that $ w \mu(w') + \mu(w)w' - \mu(ww') = 0.$ Therefore there is a
smooth vector field $\zeta \in \X(M)$ such that $ \mu(w) = \langle \ dw, \ \zeta \
\rangle$. We also have $ a \lambda(a') - \lambda(aa') = 0$. We have already shown
that such $\lambda$ is a 1-closed parallel 1-form on the locally flat manifold
$(M,D).$ The condition $\delta \theta(0,w),(a',0) = 0$ implies that $ w \lambda(a')
= 0$ for all $w \in C^{\infty}(M,\R)$ and $a' \in \X(M)$. Hence $\lambda = 0$. So
the 1-cocycle $\theta(a,w)$ has the form $\theta(a,w) = ([\xi,a], \ \langle \
dw,\zeta \ \rangle)$ for some fixed $(\xi, \zeta) \in J(\A) \times \A$. It is easy
to verify that every 1-chain $\theta$ which has the form
$$
(a,w) \longrightarrow ([\xi,a], \ \langle \ dw, \zeta \ \rangle )
$$
with $\xi \in J(\A)$ is an 1-cocycle in $C_1(\A \oplus W, \A \oplus W)$ We just
computed $H^0(\A \oplus W, \A \oplus W) \ {\rm and} \ H^1(\A \oplus W, \A \oplus W)$
:
$$
\begin{array}{l}
H^0(\A \oplus W, \A \oplus W) = \{0\};\\
H^1(\A \oplus W; \A \oplus W) \simeq \A/J(\A).
\end{array}
$$

If the manifold $M$ is compact then the dimension of $H^1(\A \oplus W, \A \oplus W)$
is infinite. The same conclusion holds if $(M,D)$ is geodesically complete.

\subsection{KV-Cohomology with values in $L(W,V)$}

Let $W$ and $V$ be KV-modules for the same KV-algebra $\A$. The vector space $W$
carries a trivial structure of KV-algebra whose multiplication map is the zero map.

We shall equip $\A \oplus W$ with the following multiplication map :
$$
(a,w)(a',w') = (aa', aw' + wa')
$$
Then $\A \oplus W$ becomes a KV-algebra such that the sequence $ 0 \longrightarrow W
\hookrightarrow \A \oplus W \longrightarrow \A \longmapsto 0 $ is an exact sequence
of KV-algebras. Furthermore the vector space $V$ becomes a KV-module for the
KV-algebra $\A \oplus W$ under the actions $ (a,w).v = av, \  v.(a,w) = va.$

Consider the KV-chain complex $C_{\tau}(\A \oplus W, V)$. It is easy to see that
$$
J_{\A \oplus W}(V) = J_{\A}(V).
$$
Now let us equip $C_{\tau}(\A \oplus W, V)$ with the bigraduation defined by
$$
C_{p,q}(\A \oplus W, V) = L(\A^q \otimes W^p,V).
$$
In other words,  elements $\theta \ {\rm in} \ C_{p,q}(\A \oplus W, V)$ are those
$(p + q)$-linear maps from $\A \oplus W$ to $V$ that are homogeneous of degree p
with respect to the elements in $W$, and homogeneous of degree $q$ in the arguments
in $\A$. By setting $k = p + q$ the vector space of $V$-valued k-chains of $\A
\oplus W$ is bigraded by the $C_{p,q}$ with $p + q = k.$ We have the following lemma
:

\begin{lem} The coboundary operator $\delta : C_k(\A \oplus W, V) \longrightarrow C_{k+1}(\A
\oplus W,V)$ maps $C_{p.q}(\A \oplus W, V)$ to \ $C_{p.q+1}(\A \oplus W, V).$
\end{lem}

It suffices to observe that if $\theta \in C_{p,q}(\A \oplus W,V)$ then
$\theta(\xi_1 \dots \xi_{p+q}) = 0$ if more than $p$ arguments belong to $W$ (resp
if more than $q$ arguments belong to $\A$). $\triangle$

The above lemma implies that $C_{\tau}(\A \oplus W, V)$ may be equipped with two
filtrations :
$$
\begin{array}{rl}
(f1) \quad F^P(\A \oplus W,V) = & {{\bigoplus} \atop {q,s \geq p}} C_{s,q}(\A \oplus W, V),\\
\\
(f2) \quad F_P(\A \oplus W,V) = & {{\bigoplus} \atop {q,s \leq p}} C_{s,q}(\A \oplus W, V).\\
\end{array}
$$

One has $F^{p+1}(\A \oplus W,V) \subset F^P(\A \oplus W,V)$ and $F_P(\A \oplus W,V)
\subset F_{p+1}(\A \oplus W,V).$ Moreover, $ \delta F^P(\A \oplus W,V) \subset
F^P(\A \oplus W,V)$ and $ \delta F_P(\A \oplus W,V) \subset F_P(\A \oplus W,V).$

Hence both filtrations will give rise to spectral sequences in the category of
``$\F$-projective" KV-modules for the KV-algebra $\A$. The study of those spectral
sequence is not the purpose of the present work. We shall mainly be interested in
the subcomplex
$$
\bigoplus_q C_{1,q}(\A \oplus W,V).
$$
The vector space $C_{1,q}(\A \oplus W,V)$ is nothing but the space of $q$-linear
mappings from $\A$ to the linear space $L(W,V)$. The vector space $L(W,V)$ is a
KV-module for $\A$ under the two actions
$$
\begin{array}{l}
(a \theta)(w) =  a_V(\theta(w)) - \theta (a_Ww),\\
(\theta a)(w) =  (\theta(w))a_V.
\end{array} \leqno (15')
$$
By identifying $L(W,V)$ with $C_{1,0}(\A \oplus W,V)$, we define
$$
\delta : C_{1,0}(\A \oplus W,V) \longrightarrow C_{1,1}(\A \oplus W,V)
$$
by putting :
$$
\left \{
\begin{array}{l}
\delta \theta(a,w) = -a \theta(w) + \theta(aw),\\
\delta \theta(w,a) = \theta(wa) - (\theta(w))a.\\
\end{array} \right. \leqno(16)
$$

It is easy to verify that ($\delta \circ \delta) \ \theta = 0$, for $\theta \in
C_{1,0}(\A \oplus W,V)$. On the other hand, each vector space $C_{1,q}(\A \oplus
W,V)$ is a KV-module for $\A \oplus W$. So we may extend the formula (4) to
${{\bigoplus} \atop {q \geq 0}} C_{1,q}(\A \oplus W,V)$ and obtain a chain complex
$$
... \longrightarrow C_{1,q}(\A \oplus W,V) \stackrel{\delta}{\longrightarrow}
C_{1,q+1}(\A \oplus W,V) \stackrel{\delta}{\longrightarrow} ... \leqno (16')
$$
\begin{df} We set $E_1^{1,q}(\A \oplus W,V)$ for the cohomology space of (16') at the level
$C_{1,q}(\A \oplus W,V)$.
\end{df}

\underline{Remark}. There exists a linear map : $E_{1,}^{1,q}(\A \oplus W,V)
\longrightarrow H^q(\A, L(W,V))$. Indeed considering the actions (15') we know that
the total chain complex of the $\A \oplus W$-KV-module $L(W,V)$ is

$$
C_{\tau}(\A \oplus W,V) = J(L(W,V) \oplus \sum_{q \rangle 0} C_q(\A \oplus W, L(W,V)).
$$

The cohomology spaces of the complex ${\oplus}_{q > 0} C_{1,q}(\A \oplus W,V)$ are
related to the spectral sequence which is associated to the filtration of
$C_{\tau}(\A \oplus W,V)$ by the $F^P(\A \oplus W,V)$.

\subsection{Consistency}

Let $W$ and $V$ be KV-modules of a KV-algebra $\A$, and $\phi : W \longrightarrow V$
be a KV-morphism. Then $\phi$ is an $\F$-linear map such that
$$
\phi(aw) = a \phi(w) \ {\rm and} \ \phi(wa) = \phi(w)a \leqno(17)
$$
for all $w \in W$ and $a \in \A$. If $f \in C_q(\A,W)$, then $\phi^*(f) \in
C_q(\A,V)$ is defined by the classical conditions : $\phi^*(f) = \phi \circ f.$
Therefore we see that

$$
\left \{
\begin{array}{l}
[ a(\phi \circ f) ] (b) = [\phi (af)] (b) \cr [(\phi \circ f) a]
(b) = [\phi (fa)] (b)
\end{array} \right.  \leqno (18)
$$

Relations (18) show that for every $f \in C_q(\A,W)$ we have $ \delta [\phi \circ f]
= \phi \circ (\delta f) $. Hence $\phi$ induces canonically the linear map $
\tilde{\phi} : H^q(\A,W) \longrightarrow H^q(\A,V).$

\underline{$\F$-splittable exact sequences}. We recall that $\F$ is a fixed
commutative field of characteristic zero. Let $W, V, T$ be three $\F$-vector spaces
which are also KV-modules for a fixed KV-algebra $\A$. Suppose that we have a short
exact sequence of KV-modules.

$$
(s) \quad 0 \longrightarrow V \longrightarrow T \longrightarrow W \longrightarrow 0.
$$

We say that (s) is $\F$-splittable if the subspace V is a direct summand in T.

If (s) is splittable then the following sequence
$$
(s') \quad 0 \longrightarrow C_q(\A,V) \longrightarrow C_q(\A, T) \longrightarrow C_q(\A, W)
\longrightarrow 0
$$
is  exact and $\F$-splittable also. From (19) we obtain the exact sequence of
cochain complexes
$$
0 \longrightarrow C_{\tau}(\A, V) \stackrel{i}{\longrightarrow} C_{\tau}(\A,T)
\stackrel{\mu}{\longrightarrow} C_{\tau}(\A,W) \longrightarrow 0
$$
which is $\F$-splittable. Therefore (s) gises rise to the long exact sequence
$$
... \stackrel{\delta}{\longrightarrow} H^q(\A,W) \stackrel{i}{\longrightarrow}
H^q(\A,T) \stackrel{\mu}{\longrightarrow} H^q(\A,W)
\stackrel{\delta}{\longrightarrow} H_{q+1}(\A, V) \stackrel{i}{\longrightarrow} ...
.
$$
Now suppose that the following map $\phi$ is a morphism of $\F$-splittable exact
sequence of KV-modules
$$
\begin{array}{llll}
0 \longrightarrow & V \stackrel{i}{\longrightarrow} & T \stackrel{P}{\longrightarrow} & W
\longrightarrow 0\\
 & \downarrow \phi & \downarrow \phi & \downarrow \bar{\phi}\\
0 \longrightarrow & V' \stackrel{i'}{\longrightarrow} & T' \stackrel{P'}{\longrightarrow} & W'
\longrightarrow 0
\end{array}
$$
Then $\phi$ will induce a morphism of long exact sequence
$$
\begin{array}{llllllll}
\stackrel{\delta}{\longrightarrow} & H^q(\A,V) & \stackrel{i}{\longrightarrow} &
H^q(\A,T) & \stackrel{P}{\longrightarrow} & H^q( \A,W) &
\stackrel{\delta}{\longrightarrow} & H_{q+1}(\A,V)
\longrightarrow\\
 & \downarrow \tilde{\phi} & &\downarrow \tilde{\phi} & & \downarrow \bar{\phi} & & \downarrow
\tilde{\phi}\\
\stackrel{\delta'}{\longrightarrow} & H^q(\A,V') & \stackrel{i'}{\longrightarrow} &
H^q(\A,T') & \stackrel{P'}{\longrightarrow} & H^q(\A,W') &
\stackrel{\delta'}{\longrightarrow} & H_{q+1}(\A,V') \longrightarrow .
\end{array}
$$

Those properties show the consistency of the KV-cohomology theory which is defined
following (4).

\section{The significane of some KV-cohomology spaces}

In this section we shall be concerned with the status of some cohomology spaces such
as $H^0(\ ), H^1(\ )$ and $H^2(\ )$. Our main aim will be to prove that those spaces
play an essential role in some important questions.

\subsection{Extensions of KV-algebras}

To begin, it is useful to recall that only some types of extension of Lie algebras
are closely related to the Hochschild cohomology theory of Lie algebra. If
$$
(e) \quad 0 \longrightarrow \B \longrightarrow \G \longrightarrow \A \longrightarrow 0
$$
is an exact sequence of Lie algebras, then $\B$ will not necessarily be an
$\A$-module, so that the equivalence class of (e) is not related to some cohomology
class of $\A$ with values in $\B$.

Now let us consider an exact sequence of KV-modules for $\A$

$$
(\tilde{s}) \quad 0 \longrightarrow W \stackrel{i}{\longrightarrow} T
\stackrel{\mu}{\longrightarrow} \A \longrightarrow 0.
$$

That is to say that $\A$ is a KV-algebra, $W$ and $T$ are KV-modules for $\A$. Let
us suppose the sequence (s) to be splittable. Then we shall identify $T$ with $W
\oplus \A$ by using a section $\A \longrightarrow T$.

From now on we make the assumption that $T$ admits a structure of KV-algebra and $W$
is a bilateral KV-ideal in $T$ the multiplication in which is the zero map.
Therefore the multiplication in $T$ is given as follows
$$
(a + w)(a' + w')  = aa' + aw' + wa'+ \omega(a,a'),
$$
where $\omega : \A \times \A \longrightarrow W$ is bilinear.

Because of the property telling that for all $a + w,a' + w' \ {\rm and}\  a'' + w''$
in $\A \oplus W$ one has
$$
(a + w, a' + w', a'' + w'') = (a' + w', a + w, a''+ w''),
$$
one easily deduces that $\omega$ is a 2-cocycle in $C_2(\A, W)$, i.e.
$$
\delta \omega(a,a',a'') = 0 \leqno (20)
$$
for all $a, a' \ {\rm and} \ a''$ in $\A$. If one identifies $T$ with $W \oplus \A$
using another section $\A \longrightarrow T$, then the induced 2-cocycle $\omega_0$
will be related to the precedent $\omega$ by
$$
\omega_0 = \omega + \delta \psi , \leqno (21)
$$
where $\psi \in C_1(\A,W)$. Thus the cohomology class $[\omega] \in H^2(\A,W)$ does
not depend on the choice of a section $\A \longrightarrow T$.

If one starts with the exact sequence of KV-algebras
$$
(\tilde{s}) \quad 0 \longrightarrow W \longrightarrow T \longrightarrow \A
\longrightarrow 0
$$
such that the multiplication in $W$ is a non-zero map, then $(\tilde{s})$ will not
necessarily be an exact sequence of $\A$-KV-modules, even if $(\tilde{s})$ is
splittable. So that $(\tilde{s})$ is not related to any 2-cohomology class of $\A$.
The similar situation is well known for non abelian extensions of Lie algebras [BN].

Another reason why we are interested in splittable exact sequences of KV-modules is
that many important examples that come from the differential geometry involve
infinite dimensional vector spaces where a vector subspace is not necessarily a
summand factor.

\subsection{$\F$-projective KV-modules}

The notion to be introduced below is motivated by the last remarks we just made. Let
$W$ be a KV-module for a KV-algebra $\A$. One says that W is {\bf{$\F$-projective}}
if every exact sequence of $\F$-vector spaces

$$
0 \longrightarrow V \longrightarrow T \longrightarrow W \longrightarrow 0
$$

is splittable. A KV-algebra $\A$ is {\bf{$\F$-projective}} if it is $\F$-projective as a
KV-module of itself.

Suppose that $\A$ is a $\F$-projective KV-algebra. Then for every KV-module $W$
equipped with the trivial multiplication map $(w,w') \longmapsto 0$, every exact
sequence of KV-modules $ 0 \longrightarrow W \longrightarrow T \longrightarrow \A
\longrightarrow 0$ is an exact sequence of KV-algebras with the multiplication
$$
(w,a)(w',\dot{a}) = (wa' + aw', aa') .
$$
Mutatis mutandis the same property holds in the case where $\G$ is an associative
algebra and $\rho : \G \longmapsto End(W)$ is a homomorphism of associative
algebras.

In the case of KV-algebras and KV-modules the analogue to the bijective map : $H^1(\G, L(W,V))
\longrightarrow$ \{equivalence classes of extensions of \ W by V\} fails. In order to bring
under control the classification problem of extensions
$$
0 \longrightarrow V \longrightarrow T \longrightarrow W \longrightarrow 0,
$$
we will study  the so-called first step of a spectral sequence which is defined by
the filtration $(f_2)$, (see the paragraph 4).

Let $W$ be a $\F$-projective KV-module for $\A$. One shall consider $W$ as a trivial
KV-algebra, so that we get the semi-direct product of KV-algebras $\A \times W$ with
the multiplication
$$
(a,w).(a',w') = (aa',aw' + wa')
$$
for $(a,w)$ and $(a',w')$ in $\A \times W \simeq \A \oplus W$.

Suppose that there is given an exact sequence of KV-modules
$$
0 \longrightarrow V \longrightarrow T \longrightarrow W \longrightarrow 0,
$$
where $W$ is a $\F$-projective module. We may identify $T$ with $V \oplus W$.

Now one considers $V$ as a KV-module of $\A \times W$ by putting $ (a,w)v = av,
v(a,w) = va.$

We have already defined the chain complex $C_{\tau}(\A \times W,V) = \sum_{q \geq 0}
C_{1,q}(\A \times W,V)$, i.e.
$$
... \longrightarrow C_{1,q-1}(\A \times W,V) \stackrel{\delta}{\longrightarrow}
C_{1,q}(\A \times W,V) \stackrel{\delta}{\longrightarrow} C_{1,q+1}(\A \times W,V)
\longrightarrow ... .
$$

We are now in position to classify the extensions of $\A$ by the trivial KV-algebra
$W$.

Indeed, let $\sigma: \A \longmapsto T$ be a section of
$$
0 \longrightarrow W \longrightarrow T \longrightarrow \A \longrightarrow 0.
$$
Consider the following 2-chain in $C_2(\A,W)$ :

$$
(a,a') \longrightarrow \omega(a,a') = \sigma(a) \sigma(a') - \sigma(aa').
$$

The cochain $\omega$ is a 2-cocycle. Thus we can define the following multiplication

$$
(w,a).(w',a') = (aw' + wa' + \omega (a,a'), aa')
$$

We see that the cohomology class $[\omega] \in H^2(\A,W)$ doesn't depend on the
choice of the section $\sigma$. Thus we can state the following classification
theorem.

\begin{Th} Let $W$ be a trivial KV-algebra that is a KV-module for $\A$. Then there is an one to one
map from $H^2(\A,W)$ onto the set of equivalent classes of splittable extensions of
$\A$ by $W$.
\end{Th}

\begin{cor} Under the assumptions of Theorem 5.1, if the KV-algebra $\A$ is $\F$-projective,
then there is a bijective correspondence between $H^2(\A,W)$ and the set of
equivalent classes of extensions of $\A$ by W.
\end{cor}

\subsection{Extensions of KV-modules}

Let $\G$ be a Lie-algebra and let $V$ and $W$ $\G$-modules, then each  splittable
exact sequence of $\G$-modules
$$
(s) 0 \longrightarrow V \longrightarrow T \longrightarrow W \longrightarrow 0
$$
will give rise to a cohomology class in $H^1(\G, L(W,V))$ which determines (s) up to
equivalence of extensions classes.

In the case of KV-modules $V,W$ of a KV-algebra $\A$, we will show that the
cohomology space responsible for equivalent classes of extension is not $H^1(\A,
L(W,V))$, but rather the space $E^{1,1}_1(\A \oplus W, V)$ at the level $C_{1,1}(\A
\times W, V)$.

Before proceeding let us fixe some notations.

Let $f \in C_{1,1}(\A \oplus W,V)$, then $\delta f$ belongs to $C_{1,2}(\A \oplus W, V)$.
For every $(a,w) \in \A \times W$ we set

$$
\begin{array}{ll}
\theta (a,w) = f(a,w),\\
\psi(a,w) = f(w,a).
\end{array} \leqno (22)
$$

Since $\delta f \in C_{1,2}(\A \oplus W,V)$ for $a, b \in \A$ and $w \in W$, one has

$$
\begin{array}{ll}
\delta f(a,b,w) & = \delta \theta(a,b,w)\\
& = -a \theta (b,w) + \theta(ab,w) + \theta(b,aw)\\
& +b \theta(a,w) - \theta(ba,w) - \theta(a,bw) ;\\
\\
\delta f(a,w,b) & = -af(w,b) + f(aw,b) + f(w,ab) - f(w,a)b\\
& +wf(a,b) - f(wa,b) - f(a,wb) + f(a,w)b.
\end{array}
$$

So that by the vertu of (22) we see that $\theta \ {\rm and} \ \psi$ are related by

$$
\left \{
\begin{array}{ll}
\delta f(a,w,b) & = -a \psi(b,w) + \psi(b,aw) + \psi(ab,w) - (\psi(a,w))b\\
& - \psi(b,wa) - \theta(a,wb) + (\theta(a,w))b.\\
\\
\delta f(a,b,w) & = \delta \theta(a,b,w).
\end{array} \right. \leqno (23)
$$

So that when $f \in C_{1,2}(\A \times W,V)$ is a cocycle we get

$$
\left \{
\begin{array}{l}
\delta \theta (a,b,w) = 0\\
a \psi(b,w) + \psi(b,aw) + \psi(ab,w) - (\psi(a,w))b - \psi(b,wa) = \theta(a,wb) -
(\theta(a,w))b.
\end{array} \right. \leqno (23')
$$

Let us consider an exact sequence of KV-modules

$$
(s) \quad 0\rightarrow V \longrightarrow T \longrightarrow W \longrightarrow 0.
$$

Fix a section $\sigma : W \longrightarrow T. Then \sigma$ will define the following
(1,1)-cochain in $C_{1,1}(\A \times W,V)$:

$$
\begin{array}{l}
f_{\sigma}(a,w) = a \sigma (w) - \sigma (aw)\\
f_{\sigma}(w,a) = \sigma (w)a - \sigma (wa)
\end{array} \leqno (24)
$$

\begin{lem} $\delta f_{\sigma} = 0$.
\end{lem}

Proof. Since $\delta f_{\sigma}$ is a (1,2)-chain of $\A \times W$ with values in
$V$, for any $a \ {\rm and} \ b$ in $\A$ and $w \in W$ we have :

$$
\begin{array}{ll}
\delta f_{\sigma}(a,b,w) & = -a(b \sigma (w) - \sigma (bx)) + (ab) \sigma (w) - \sigma((ab)w)\\
 & + b \sigma (aw) - \sigma (b(aw)) - b(a \sigma(w)) + b \sigma (aw)\\
 & -(ba) \sigma(w) + \sigma((ba)w) - a\sigma(bw) + \sigma(a(bw))\\
 & = (b,a,\sigma(w)) - (a,b,\sigma(w)) + \sigma((a,b,w) - (b,a,w))\\
 & = 0
\end{array} \leqno (i)
$$

$$
\begin{array}{ll}
\delta f_{\sigma}(a,w,b) & = -a(\sigma(w)b) - a \sigma(wb) + \sigma(aw)b - \sigma((aw)b)\\
 & + \sigma(w)ab - \sigma(w(ab)) - (\sigma(w)a)b + \sigma(wa)b\\
 & - \sigma(wa)b + \sigma((wa)b) - a\sigma(wb) + \sigma(a(wb))\\
 & + (a \sigma(w))b - (\sigma(aw))b \\
 & = (a, \sigma(w),b) - (\sigma(w),a,b) + \sigma((a,w,b) - (w,a,b))\\
 & = 0
\end{array} \leqno (ii)
$$

Lemma 5.1 is proved. $\triangle$

Now
$$
\begin{array}{ll}
\theta(a,w) = f_{\sigma}(a,w) \ {\rm and}\\
\psi(a,w) = f_{\sigma}(w,a)
\end{array} \leqno (25)
$$
are related as in (24). It follows from Lemma 5.1. that the exact sequence
$$
0 \longrightarrow V \longrightarrow T \longrightarrow W \longrightarrow 0
$$
gives rise to a class $[f_{\sigma}] \in E^{1,1}_1(\A \times W, V)$, where
$f_{\sigma}$ is defined as in (24).

It is easy to verify that if $\sigma' : W \longrightarrow T$ is another section of
the exact sequence (s) then $[f_{\sigma'}] = [f_{\sigma}]$, so that $[f_{\sigma}]$
depends only on the equivalence class of the extension of $W$ by the KV-module $V$.
Of course, given a class $[f] \in E^{1,1}_1$, one can construct an extension of $W$
by $V$,
$$
0 \rightarrow V \longrightarrow T \longrightarrow W \longrightarrow 0 ,
$$
by setting
$$
\begin{array}{l}
\theta (a,w) = f(a,w),\\
\psi(a,w) = f(w,a),\\
a(v,w) = (av + \theta (a,w), aw),\\
(v,w)a = (va + \psi(a,w),wa).
\end{array}
$$
for $f \in [f] \ {\rm and} \ (a,w) \in \A \times W, ((v,w) \in T = V \oplus W)$.

We can now state the classification theorem of extensions of KV-modules, which is
analogous to the ones for associative algebras and for Lie algebras.

\begin{Th} Let $W$ and $V$ be KV-modules of a KV-algebra $\A$. Then there
is a bijective correspondence between $E^{1,1}_1(\A \times W,V)$ and the set
$Ext(W,V)$ of classes of splittable exact sequences of KV-modules $0 \longrightarrow
V \longrightarrow T \longrightarrow W \rightarrow 0.$
\end{Th}

\begin{cor} Under the notations in Theorem 5.2, if $W$ is $\F$-projective, then $ E_1^{1,1}(\A
\times W,V)$ corresponds bijectively to the set $Ext(W,V)$ of classes of exact
sequences of KV-modules $0 \longrightarrow V \longrightarrow T \longrightarrow W
\rightarrow 0.$
\end{cor}

To conclude, observe that $E^{1,1}_1(\A \times W,V)$ plays in the theory of
KV-modules the same role as $H^1(\G, L(W,V))$ in the theory of modules for
associative algebra or Lie algebras.

\section{Deformations of algebraic structures}

In the study of algebraic structures on a given vector space, the Hochschild
cohomology theory is involved either as obstruction to rigidity of some fixed
structure, or as a tool to compare Zariski neighbourhood of a fixed structure and
its orbit under the action of some group. The cases of modules of associative and
Lie algebras are well studied [GM], [NA], [KJL2]. We intend to raise the analogous
problems for KV-algebras and KV-modules. To begin with, let us fix some notations.

Let $\A$ be a KV-algebra. We denote by $|\A|$ the underlying vector space of $\A$.
Let $KV(|\A|)$ be the set of KV-algebra structures on $|\A|$. It is an algebraic
variety (eventually singular). If $\mu$ is a given element in $KV(|\A|)$ we shall
set $\mu(a,b) = ab$ when there is no risk of confusion.

Of course $KV(|\A|)$ is a subset of the space $B(|\A|, |\A|)$ of bilinear maps from
$\A$ to itself. For $\mu$ and $\nu$ in $B(|\A|,|\A|)$ we denote by $d_{\mu}\nu$ for
the three-linear map defined as follows:

$$
\begin{array}{ll}
d_{\mu}\nu(a,b,c) = & - \mu(a,\nu(b,c)) + \nu(\mu(a,b),c) + \nu(b,\mu(a,c)) - \mu(\nu(b,a),c)\\
& + \mu(b,\nu(a,c)) - \nu(\mu(b,a),c) - \nu(a,\mu(b,c)) + \mu(\nu(a,b),c).
\end{array} \leqno (26)
$$

To avoid confusion with the Gerstenhaber-Nijenhuis bracket, we don't write $[\mu,
\nu].$

\subsection{Deformations of KV-algebras}

Let $\A$ be a KV-algebra. The multiplication in $\A$ will be denoted by $(a,b)
\longmapsto ab.$ Suppose that there is a one parameter family of KV-algebra
structures on $|\A|$ : $ \mu(t) = \sum^{\infty}_{j=0} \mu_jt^j$  with
$$
\begin{array}{l}
\mu_o(a,a') = aa',\\
\mu_j \in B(|\A|,|\A|).
\end{array}
$$
Henceforth we shall set $(a,a',a'')_t$ for the $\mu(t)$-associator. Because $(a,a')
\longmapsto \mu_t(a,a')$ defines a KV-algebra structure in $|\A|$ we have
$$
(a,b,c)_t - (b,a,c)_t = 0 \leqno (27)
$$
for all $a, b$ and $c$ in $|\A|.$ Using notation (26), identity (27) is equivalent
to the following
$$
\begin{array}{l}
\delta \mu_1 = 0\\
\delta \mu_k = \sum_{{i+j=k} \atop{i \rangle 0,j \rangle 0}} \delta_i \mu_j
\end{array} \leqno (27')
$$
where $\delta_i \mu_j = \delta_{\mu i} \mu_j$. We identify the Zariski tangent space
of $KV(|\A|)$ at $\A \in KV(|\A|)$ (or at $\mu_o \in KV(|\A|$)) with the subspace of
2-cocycles $Z_2(\A,|\A|) = Z_2(\A,\A) \subset C_2(\A,\A)$. So the family $\mu_t$
gives rise to the cohomology class $[\mu_1] \in H^2(\A,\A)$.

Consider an one-parameter family of linear isomorphisms from $|\A|$ to itself. One
supposes $\phi_t$ to have the power series form
$$
\phi_t(a) = a + t \theta_1(a) \dots + t^j \theta_j(a),
$$
where the $\theta_j$ are linear endomorphisms of $|\A|$. Let us set
$$
\mu_t(a,a') = \phi_t(\phi^{-1}_t(a) \phi^{-1}_t(a')) .
$$
So we define a family of KV-algebra structures $\A(t) = (|\A|, \mu_t)$ which are
isomorphic to $\A$. We expand $\mu_t(a,a')$ as
$$
\mu_t(a,a') = aa' + \sum_{k \rangle 0} t^k \mu_k(a,a')
$$
with $\mu_k \in B(|\A|,|\A|)$. Therefore we use the relation (27) to see that
$$
\mu_1(a,a') = \delta \theta_1(a,a').
$$
So we see that $H^2(\A,\A)$ may be regarded as the Zariski tangent space of
$KV(|\A|)$ at $\A$. In particular we get

\begin{Th} (the rigidity theorem). If
$H^2(\A,\A)$ is trivial then $\A$ is formally rigid.
\end{Th}

The space $H^2(\A, \A)$ can be viewed as the set of non trivial infinitesimal
deformations of $\A$. In some special situation, $H^2(\A,\A) \not= 0$ will imply
that $\A$ admits non trivial deformations (see Part II).

\part{Differential geometry}

In this Part II, we intend to raise a few problems that motivate the intrinsic
KV-cohomology theory. One may observe the following facts.

1) Associative algebras admit another cohomology theory that is different from the
Hochschild cohomology theory. Indeed, if $\A$ is an associative algebra then it is a
KV-algebra with $J(\A) = \A$. Every cocycle of degree 2 of the associative algebra
$w \in C^2_{ass}(\A,\A)$ is also a 2-cocycle for the KV-algebra $\A$. The inverse is
not true, so one gets a homomorphism
$$
H^2(\A,\A)_{ass} \longrightarrow H^2(\A,\A)_{KV}.
$$
which is not surjective.

2) The KV-cohomology theory developed in Part I is different from the one initiated
by A. Nijenhuis, [NA]. Nijenhuis' theory is nothing but the Hochschild cohomology
theory for Lie algebras with coefficients in spaces of linear maps.

Let us return to the differential geometry.

\section{KV-Cohomology of locally flat manifolds}

Let $M$ be a smooth manifold equipped with a torsion free-linear flat connection
$D$. Let $\X(M)$ be the vector space of smooth vector fields on $M$ and $TM$ be the
tangent vector bundle.  The connection $D$ give rise to the KV-algebra
$$
\A = (\X(M), D), ab = D_ab \ ; a, b \ {\rm in} \ \X(M)).
$$
Let $T(M)$ denote the space of tensors on $M$. Then $T(M)$ is an infinite
dimensional bi-graded vector space :
$$
T(M) = \sum_{p,q} T^{p,q}(M) ,
$$
where $T^{p,q}(M)$ is the vector space of smooth sections of the vector bundle
$T^{p,q}M = (\stackrel{p}{\otimes} TM) \otimes (\stackrel{q}{\otimes}T^*M$); $T^*M$
is the cotangent bundle of $M$. The vector space $T(M)$ is a KV-module of the
KV-algebra $\A = (\X(M), D)$. To see that, it is sufficient to specify the actions
of $\A$ on $T^{1,0}(M) \ \hbox{\rm and} \ T^{0,1}(M)$. If $a \in \A, X \in
T^{1,0}(M) \ {\rm and}\ \theta \in T^{0,1}(M)$ then we shall set

$$
 aX = D_aX, \ Xa = D_Xa, \ a \theta = D_a \theta \ , \ \theta a = 0.
$$

The differential 1-form $D_a \theta$ is defined by the formula
$$
< (D_a \theta),X > \ = \ < d(\theta(X)),a > - < \theta,aX >.
$$
Of course $T^{p,o}(M)$ is defined for $p > 0$ while $T^{o,q}$ is defined for non
negative integers under the convention that $ T^{0,0}(M) = C^{\infty}(M, \R). $ We
can define the KV-chain complex
$$
C(\A, T(M)) = \bigoplus_\rho C_{\rho}(\A,T(M))
$$
whose boundary operator
$
\delta : C_{\rho}(\A,T(M)) \longrightarrow C_{\rho +1}(\A,T(M))
$
is defined by the formula (4), viz, taking $f \in C_{\rho}(\A,T(M))$ we have
$$
\delta f(a_1 \dots a_{\rho +1}) = \sum_{i \leq \rho}(-1)^i \{(a_if)(\dots \hat{a}_i \dots
a_{\rho +1}) + (e_{\rho} (a_i).(fa_{\rho +1})(\dots \hat{a}_i \dots a_{\rho})\}.
$$
Following Lemma 4.1 we get the following chain complex
$$
... \stackrel{\delta}{\longrightarrow} C_{p -1}(\A, T(M)) \stackrel
{\delta}{\longrightarrow} C_{p}(\A, T(M)) \stackrel{\delta}{\longrightarrow} C_{p
+1}(\A,T(M)) \longrightarrow .
$$
The $p^{th}$ cohomology space
$$
\frac{ker \ \{\delta : C_p(\A,T(M)) \longrightarrow C_{p +1}(\A,T(M))\}} {\delta
C_{p -1}(\A,T(M))}
$$
is denoted by $H^p(\A,T(M))$. It is likely that the spaces $H^p(\A, T(M))$ are very
big. It is also to be noticed that in general cochains in $C_{\rho}(\A,T(M))$ are
not tensors (viz $C^{\infty}(M, \R)$-multilinear), but they are tensor-valued
$\R$-multilinear maps. One may also observe that $C(\A,T(M))$ admits a triple
graduation :

$$
C(\A,T(M)) = \sum C^{p,q}_{\rho}(\A,T(M)).
$$

Indeed the actions of $\A$ in $T(M) = {{\oplus}_{p,q}} T^{p,q}(M)$ has degree zero,
hence every subspace $T^{p,q}(M)$ is a KV-module for $\A$. Therefore we can write

$$
C_{\rho}^{p,q}(\A,T(M)) = C_{\rho}(\A,T^{p,q}(M)).
$$

Then the complex $C(\A,T(M))$becomes a ``bi-graded" chain complex in the sense that
it is a direct sum of a two-parameter family of subcomplexes $C^{p,q}(\A,T(M)) =
{{\oplus}_{\rho}} C_{\rho}^{p,q}(\A,T(M))$.

If we set $W = T^{0,0}(M) = C^{\infty}(M, \R)$ we get $C(\A,W)$ as a subcomplex of
$C(\A,T(M))$. For every non negative integer $\rho$ the space of tensorial $\rho$-forms $\tau
^p(M, \R)$ is a subspace of $C_{\rho}(\A,W)$, it consists of chains $f \in C_{\rho}(\A,W)$ which
are tensorial. It is easy to verify that the subspace $\tau (M,\R)$ of tensorial forms is
$\delta$-stable. Hence $\tau(M,\R)$ is a subcomplex of the complex $C(\A,W)$.

\begin{df} The subcomplex $\tau(M,\R) \ {\rm of} \ C(\A,W)$ is called the KV-complex of the
locally flat manifold $(M,D)$.
\end{df}

\underline{Remark}. In [KJL1] Jean-Louis Koszul defines the $D$-complex of
$T*M$-valued forms. Indeed one can use the locally flat connection $D$ to get an
exterior differential operator in $T^*M$.

To do so one regards D as a linear representation of the Lie algebra of smooth
vector fields $\{\X(M)$, Lie bracket\} in $\X(M)$. That gives rise to the dual
representation in $T^{0,*}(M)$. In particular the Hochschild complex that one gets
will give rise to the $TM$-valued (resp $T^*M$-valued) cohomology. What has been
done in [NA] is the algebraic version of that procedure. What we are going to do is
quite different.

We call the subcomplex $(\tau(M,\R),\delta)$ the KV-complex of $(M,D)$. So the
complex $C(\A, T(M))$ contains as subcomplex the two complexes $\tau(M,\R) \ {\rm
and} \ C(\A,W)$ :

$$
\tau(M,\R) \subset C(\A,W) \subset C(\A,T(M)).
$$

The cohomology $H_{KV}(M,\R) = \oplus_{\rho} H^{\rho}(\tau(M,\R))$ is called the
KV-cohomology of $(M,D)$. It is an affine invariant of the affine manifold $(M,D)$,
and the claim, viz to be affine and invariant, holds for the cohomology space of the
complexes $(C(\A,W) \ {\rm and}\ C(\A,T(M))$.

\section{Rigid hyperbolic affine manifolds}

We just associated to every locally flat manifold $(M,D)$ the complex $C(\A,T(M))$
containing two canonical subcomplexes : $\tau(M,\R) \ {\rm and} \ C(\A, W)$. An
interesting question is to give a geometrical meaning to the cohomology spaces of
those complexes. In the following subsection we are going to give a partial answer
to this question.

\subsection{Hyperbolic affine manifolds}

\begin{df} A locally flat manifold $(M,D)$ is called hyperbolic if its universal covering
$(\tilde{M} ,\tilde{D})$ is isomorphic to a convex domain in the Euclidian space
which does not contain any straight line.
\end{df}

Let $M$ be a smooth manifold admitting a hyperbolic locally flat structure. We
denote by $\f(M)$ the set of locally flat connections on $M$ and $C(M)$ the subset
consisting of those $D_0 \in \f(M)$ such that the universal covering of $(M,D_0)$ is
isomorphic to a convex cone not containing any straight line

In [KJL2] Jean Louis Koszul has proved the following

\begin{Th} [[KJL2] Theorem 3] If M is compact then $C(M)$ is an open subset
in $\f(M)$.
\end{Th}

In [KJL2] the proof of that theorem involves locally Hessian Riemannian metric
(under the assumption that $C(M) \not= \0$ such a metric exists).

A remarkable consequence of the above theorem is the non rigidity of the $(M,D_0)$
with $D_0 \in C(M)$. Indeed J. L. Koszul proved that if $C(M)$ is non empty then for
$D_0 \in C(M)$ the hyperbolic structure $(M,D_0)$ admits non trivial deformations.
Roughly speaking every neighbourhood of $D_0$ in $\f(M)$ contains $D$ such the
locally flat manifold $(M,D)$ is not isomorphic to $(M,D_0)$.

Using the complex $(C(\A,\A)$ where $\A = (\X(M), D_0))$, the rigidity problem for
$(M,D_0)$ is related to $H^2(\A,\A)$. From this viewpoint, the compactness of $M$
may be irrelevant. We plan to highlight our remark by studying some hyperbolic
affine structure in the manifolds $\R^2$ and $A f f(\R)$ where $A ff(\R)$ is the
group of affine transformations of the straight line.

\underline{Example}

We consider $\R^2 = \{\lambda \frac{\delta}{\delta x} + \mu \frac{\delta}{\delta y} \ ,\
\lambda \in \R \ , \ \mu \in \R\}$ with the following $\R$-bilinear bracket

$$
[\frac{\delta}{\delta x} \ , \ \frac{\delta}{\delta y}] = \frac{\delta}{\delta y}.
$$
Therefore we obtain the Lie algebra $A ff(\R)$ of infinitesimal affine
transformations of $\R$; it is the Lie algebra of the Lie group $A ff(\R)$, so that
by setting
$$
D_{(\lambda,\mu)}(\lambda',\mu')  = \lambda \mu'
\frac{\delta}{\delta y} = (0, \lambda \mu') \leqno (28)
$$

one gets a left invariant linear connection on the Lie group $A ff(\R)$.

On the manifold $A ff(\R)$ we consider the tensor $S_{\alpha, \beta} \in T^{1,2}(M)$
which is defined as  follows :
$$
S_{\alpha,\beta}((\lambda,\mu)\ ,\ (\lambda',\mu')) = (\alpha \lambda \lambda'\ ,\
\beta \lambda \lambda' \ + \ \alpha(\lambda \mu' + \mu \lambda')),
$$
where $\alpha \ {\rm and} \ \beta$ are fixed real numbers, $\lambda, \mu, \lambda' \
{\rm and} \ \mu'$ are elements of $W = C^{\infty}(A ff(\R))$. We extend canonically
$D$ to get a linear connection on $A ff(\R)$, which is locally flat, so that we get
the KV-algebra $\A = (X(\A ff(\R),D)) $.

\begin{lem} For every pair of real number $(\alpha,\beta)$ with $\alpha \not= 0$ the
$\R$-bilinear map $S_{\alpha \beta} : ((\lambda,\mu),(\lambda',\mu'))
\longrightarrow (\alpha \lambda \lambda', \beta \lambda \lambda' + \alpha(\lambda
\mu' + \mu \lambda'))$ is a non trivial 2-cocycle in $C_2(\A,\A)$ such that
$\delta_{s_{\alpha \beta}} S_{\alpha \beta} = 0$.
\end{lem}

\underline{Remarks.}

1) a cocycle is called \underline{non trivial} when it is non exact.

2) the datum $\delta_{S_{\alpha \beta}} S_{\alpha \beta}$ is defined as in (26).

The proof of Lemma 8.1 is easy. We deduce from it the following consequence.

\begin{cor} The locally flat manifold $(A ff(\R),D)$ where D is defined by (28) admits non
trivial deformations.  \end{cor}

\underline{Proof}. We define the following one parameter family of linear
connections in $A ff(\R)$ :

$$
D_t = D + tS_{\alpha \beta} .
$$

By the vertu of the lemma 8.1, $D_t$ is locally flat. We deduce the one parameter
family of KV-algebras

$$
\A(t) = (\X(\A ff(\R)),D_t).
$$

Since $S_{\alpha \beta}$ is a non trivial cocycle, the family $\A(t)$ is a non
trivial deformation of $\A = \A(0)$. This ends the proof.

\underline{Remark}. One can show that $\A(t)$ is a non trivial deformation of $\A$ by
considering the $\R$-linear map

$$
D_t(\lambda_0,\mu_0) : (\lambda,\mu) \longrightarrow D_{\lambda,\mu}(\lambda_0
\mu_0) + S_{\alpha \beta}((\lambda, \mu),(\lambda_0,\mu_0)) + (\lambda , \mu).
$$

It shows that D = $D_0$ is geodesically complete. On the other hand, if $\alpha
\not= 0$ then $D_t$ is geodesically non complete if $t \not= 0$.

That proves that the $D_t$ are not isomorphic to each other.

We are going to examine more closely $D_1 = D + S_{\alpha \beta}$.

$$
D_1 : ((\lambda, \mu),(\lambda',\mu')) \longrightarrow (\alpha \lambda \lambda', \beta \lambda
\lambda' + (1 + \alpha)\lambda \mu' + \alpha \lambda' \mu).
$$

Let us attempt to understand what do geodesics of $D_1$ look like.

Consider a geodesic $t \longrightarrow (x(t), y(t))$. We have
$$
D_{(\dot{x}, \dot{y})} (\dot{x}, \dot{y}) = (0,0),
$$

where $\dot{x}$ means $\frac{dx}{dt}$. Then we get $(\dot{x},\dot{y}) = \frac{dx}{dt} \
\frac{\delta}{\delta x} + \frac{dy}{dt} \ \frac{\delta}{\delta y}$.

The geodesic is subject to  the following system of equations :

$$
\begin{array}{l}
2 \frac{d^2x}{dt^2} + \alpha (\frac{dx}{dt})^2 = 0,\\
\\
2 \frac{d^2y}{dt^2} + \beta(\frac{dx}{dt})^2 + (1 +2 \alpha) \frac{dx}{dt} \ \frac{dy}{dt} = 0.
\end{array}
$$

The first equation, i.e. $ 2 \frac{d^2x}{dt^2} + \alpha (\frac{dx}{dt})^2 = 0 $,
admits the solutions $ x(t) = \frac{2}{\alpha} \ell n \ |\frac{\alpha}{2} t + u| +
v$, where $u$ and $v$ are arbitrary constants. It becomes easy to integrate the
second equation; its solutions are
$$
y(t) = \frac{\beta}{\alpha (1 -2 \alpha)} \ \ell n \ |\frac {\alpha}{2} t + u| -
\frac{c}{1+\alpha} exp (- \ \frac{1 +\alpha}{\alpha} \ \ell n \ |\frac{\alpha}{2}t + u|) + d
$$
where $c$ and $d$ are arbitrary constants.

So, geodesics of $D_1$ are defined only on half-lines $]\gamma, \infty[$  or $] -
\infty, \gamma[$, where $\gamma$ denotes an arbitrary but finite number. That means
that $D_1$ doesn't admit any complete geodesic whenever $\alpha \not= 0$. Of course,
the existence of functions $x(t)$ and $y(t)$ depends on other conditions on $\alpha$
such as $1 + \alpha \not= 0$ and so on \dots

Our purpose is to show that $(\A ff(\R), D_1)$ is isomorphic to a convex cone not
containing any straight line. We just saw that $(\A ff(\R),D_1)$ has no complete
geodesic. So we need only to prove that it is isomorphic to a convex cone.

Recalling that $[\frac{\delta}{\delta x},\frac{\delta}{\delta y}] = \frac{\delta}{\delta y} \
{\rm in} \ A ff(\R)$, one sees that the vector field $\frac{\delta}{\delta^x}$ commutes to
$e^{-x}\frac{\delta}{\delta^y}$. They are affine vector fields.

It is easy to verify that if $1 + \alpha \not= 0 \quad (\alpha \not= 0)$, then the functions

$$
\begin{array}{l}
g(x,y) = exp((1 + \alpha)x + y),\\
h(x,y) = exp((1 + \alpha)x + 2y + y^3)
\end{array}
$$

are $D_1$-affine functions, i.e we have

$$
D^2_1g = 0, D^2_1h = 0.
$$

Let us identify the unit connected component of $\A ff(\R)$ with $\{(x,y) \in \R^2 x >
0\}$; then the map

$$
(x,y) \longrightarrow (g(x,y),h(x,y))
$$

is an affine isomorphism of $A ff(\R)_0$ onto the convex cone $\{(u,v) \in \R^2/u > 0,
v > 0\}$.\ Thus for $\alpha \not= 0 \ {\rm and} \ \alpha +1 = 0, \ (A ff(\R)_0,D_1)$ is
hyperbolic and is isomorphic to a convex cone not containing any straight line.

\begin{pro} Let (M,D) be a compact locally flat manifold and $\A = (\X(M), D)$ be its KV-algebra
. If $H^2(\A,\A) = 0$ then every smooth one-parameter family deformations of (M,D)
is trivial.
\end{pro}

\underline{Proof}. Let $(M,D_t)$ be a smooth deformation of $(M,D_0 = D)$.

Let one set $S = \frac{d}{dt} D_{t/t = 0}$ ; then $S \in C_2(\A,\A)$ is a 2-cocycle.

Since $H^2(\A,\A)$ vanishes there exists a 1-chain $\phi \in C_1(\A,\A)$ such that

$$
S(a,b) = \delta \phi(a,b).
$$

The symmetry property of S implies that $\phi$ is a derivation of the Lie algebra $(X(M), [ \ ,\
])$ \ where \ [ \ , \ ] is the Lie bracket of smooth vector fields on M; thus by the vertu of
a classical theorem by Takens there is a $\xi \in \X(M)$ such that for every $a \in \A$ one has

$$
\phi(a) = [\xi,a].
$$

The 2-boundary $S = \delta \phi$ now takes the following form

$$
S(a,b) = (a,b, \xi).
$$

We now denote by $\phi_{(t)}$ the local flow of $\xi$. It is easy to see that D and $D_t$ are
related by

$$
D_t = \phi_{(t)} D \ ;
$$

so therefore, for elements $a$ and $b$ in $\A$ we may write

$$
(D_t)a(b) = \phi_{*}(t) . \ D_{\varphi_*(-t) . a} \ \varphi_*(-t). b
$$

That ends the proof of Proposition 8.1.

\section{\underline{Complex of superorder forms [K2]}}

\subsection{\underline{Differential forms of order $\leq$ k}.}

Let M be a smooth manifold. We are going to deal with the spaces of $C^{\infty}$-tensors that
are sections of the vector bundle $ M \longleftarrow \stackrel{q}{\otimes} T^*M$. Those sections
may be regarded as multilinear $C^{\infty}(M,\R)$-valued maps defined in $T^{q,0}(M)$. So let us
recall the notion of k-order forms. Let $\theta$ be a smooth section of
$\stackrel{q}{\otimes}T^*M, x_1, \dots x_q$ be q smooth vector fields, and $k$ be a non
negative integer.

\begin{df} $[KJL]_1 \ \theta$ is of order $\leq$ k if at every point $p \in M$ the value of
$\theta(x_1, \dots , x_q)$ at p depends on the k-jets at p of the $X_i$.
\end{df}

Let $E \longrightarrow M$ be a tensor vector bundle and let $\E$ be the real vector
space of smooth sections of E; $\E$ becomes a $\A(M)$-module under the Lie
derivation $L_X, X \in \A(M); \A(M)$ is the Lie algebra $(\X(M), [ \ ,\ ]).$ Let
$C(\A(M), \E)$ be the Chevalley-Hochschild complex of $\E$-valued multilinear maps
defined in $\A(M)$. In $[KJL]_1$ J. L. Koszul observed that in many situations the
cohomology space of $C(\A (M), \E)$ contains non vanishing canonical classes.

For instance:

1) The divergence class of $(M,v)$ where $v$ is a volume form.

2) For every linear connection D on M the linear map $X \longmapsto L_XD$ is a
$T^{1,2}(M)$-valued 1-form of order $\leq$ 1 which is actually a non trivial
cocycle.

Our purpose is similar to that in $[KJL]_1$; mutatis mutandis we deal with KV-complex on locally
flat manifolds.

Let (M,D) be a locally flat manifold and let $\A$ be its KV-algebra, ($\A(M)$ is nothing but the
Lie algebra of commutators of $\A$). Let W be the left KV-module $C^{\infty}(M,\R)$ \ (for $\A$).

We already observed that the complex $C(\A,W)$ contains chains of positive order,
so that the KV-complex of $M, \tau(M,R)$ consists of chains of order $\leq 0$.

For instance :

1) Symplectic forms, (pseudo) riemannian metrics may be regarded as 2-chains of order $\leq 0$.

2) If g is a (pseudo) riemannian metric, then its Levi-civita connection $\bigtriangledown$ is
2-chain of order $\leq 1$ in $C(\A,\A)$. Actually the connection $\bigtriangledown$ is defined
by the symmetric 2-chain

$$
S_{\bigtriangledown}(a,a') = \bigtriangledown^{a'}_a - aa'.
$$

Of course the boundary of the identity map $a \longmapsto a \ {\rm  in} \ \A \ {\rm is}
\ \delta_{1_{\A}}(a,a') = - aa'$. So that

$$
\delta \bigtriangledown (a,a',a'') = \delta S_{\bigtriangledown}(a,a',a'').
$$

So $\bigtriangledown$ is $\delta$-closed if and only if $S_{\bigtriangledown}$ is
$\delta$-closed. When $\delta \bigtriangledown = 0$ the curvature tensor R$\bigtriangledown$
of $\bigtriangledown$ is given by the formula

$$
R_{\bigtriangledown}(X,Y) = [S_{\bigtriangledown}(X, -), S{\bigtriangledown}(Y, -)] \leqno (*)
$$

where $S_{\bigtriangledown}(X, -) : Z \longrightarrow S_{\bigtriangledown}(X, Z)$. The formula
(*) may be interesting in some situations. For instance it may make easy the study of other
riemannian invariants which are related to the curvature tensor.

Of course, for any non negative integer k, the boundary of chains of order $\leq k$ are chains
of order $\leq k +1$.

Before proceeding further we will give some interesting examples.

1) Let $\omega$ be closed differential 2-form. We regard it as a 2-chain of order $\leq 0$ in the
KV-complex $C(\A,W)$. Then we get

$$
\delta \omega(a,b,c) = (c.\omega)(a,b).
$$

In other words we get $\delta \omega(a,b,c) = c(\omega(a,b)) - \omega(ca,b) - \omega(a,cb)$.

so that $\delta \omega = 0$ iff $\omega$ is parallel 2-form w.r.t the linear connection D.

Let us assume $\omega$ \ to be $\delta$-exact, then there is a 1-chain $\theta$ in $C_1(\A,W)$
such that

$$
\omega = \delta \theta.
$$

So, for X and Y in $\A$ we get

$$
\omega (X,Y) = -X \theta (Y) + \theta (XY)
$$

therefore we may conclude that
$$
-\omega (X,Y) = \frac{1}{2}(X \theta (Y) - Y \theta (X) + \theta ([X,Y])).
$$

So, a parallel closed 2-form  is $\delta$-exact if and only if it is ``de Rham" exact.

Here is an example where there is a non $\delta$-exact parallel closed 2-form :

2) We consider the direct product $M = H_3 \times \R$ where $H_3$ is the 3-Heiserberg group.
Let $\h_3$ be the Lie algebra of $H_3$, then we fix a base of $\h_3 \ e_1, e_2, e_3$ such that
$[e_1, e_2] = e_3$. So that M becomes a Lie group with the Lie algebra $\M = \h_3 \oplus \R$.
If we fix $\{e_1,e_2,e_3,e_4\}$ be a basis of $\h_3 \oplus \R$ and
$\{\epsilon_1,\epsilon_2,\epsilon_3,\epsilon_4\}$ its dual basis then $ \omega = \epsilon_1
\wedge \epsilon_4 + \epsilon_3 \wedge \epsilon_2$ is a left invariant closed 2-form on M. The
form $\omega$ is a sympletic form. Consider the 2-dimensional subgroups $N_1$ and $N_2$
corresponding to the subalgebras $\n_1 = span(e_1,e_3) \ {\rm and} \ \n_2 = span(e_2,e_4)$
respectively; both of them are Lagrangian submanifold of M. So there exist two lagrangian
foliations in $(M,\omega)$ which are transverse everywhere and which are left invariant as well.

Let $\Gamma$ be a lattice in M, (such a $\Gamma$ exists), the compact manifold $\Gamma
\setminus M$ carries a sympletic form and two lagrangian foliations that are transverse to each
other everywhere. There exists a unique torsion free linear connection D in $\Gamma \setminus M$
which preserves the two lagrangian foliations and such that $D\omega = 0, [NB]_4$

Now let (x,y,z,t) be the euclidian coordinates of $\M$ in the basis $(e_1,e_2,e_3,e_4)$. The
(x,t,y,z) give local coordinate functions on $\Gamma \setminus M$ which are Darboux
coordinate for $(\Gamma \setminus M,\omega)$.

Because those coordinate functions are caracteristic functions for the two
lagrangian foliations just defined, a theorem by Hess implies that the curvature
tensor of D vanishes. So $(\Gamma \setminus M,D)$ is a locally flat manifold and
$\omega$ will define a $\delta$-closed 2-cocycle in the KV-complex of $(\Gamma
\setminus M,D)$.

Because the manifold $\Gamma \setminus M$ is closed, $\omega$ can not be ``de
Rham"-exact therefore $\omega$ is not $\delta$-exact as well. So that we get the
analogue to to Koszul's remark.

3) We just saw that general locally flat manifolds (M,D) may carry D-parallel
cocycles that are non exact. In contrast with that general fact we get the following
vanishing theorem :

\begin{pro} Let (M,D) be a hyperbolic locally flat manifold. If the universal covering
$(\tilde{M}, \tilde{D})$ \ of $(M,D)$ is isomorphic to a convex cone then every D-parallel
2-chain in $C_2(\A,W)$ is $\delta-exact$.
\end{pro}

\underline{Proof}. By the results of $[KJL]_2$ (see lemme 3 ibidem) there exists $H \in \A$ such
that $aH =a$ for all $a \in \A$.

If $g$ is a 2-chain which is D-parallel, we define the 1-chain $\theta$ by setting

$$
\theta (a) = g(H,a)
$$

Since $Dg = 0$ we see that for $a$ and $b$ in $\A$ one has

$$
\begin{array}{ll}
0 & = ag(H,b) - g(aH,b) - g(H,ab)\\
 & = - \delta \theta (a,b) - g(a,b)
\end{array}
$$

That ends the proof of Proposition 9.1.

4) An other interesting instance is the following. Let $(M, D)$ be a locally flat manifold
whose the universal covering is denoted by $(\tilde{M}, \tilde{D})$. Let $\Gamma$ be the
fundamental group of M : the linear holonomy of $(M, D)$ is denoted by $\ell(\Gamma), [FGH],
[CY]$

If $\ell(\Gamma)$ is unimodular (viz $\ell(\Gamma)$ is included in the special
linear group,) then $\tilde{M}$ carries a volume form $v$ which is
$\tilde{(D)}$-parallel ; hence $\tilde{(v)}$ is a cocycle of the complex
$C(\tilde{\A}, \tilde{W})$ where $\tilde{(\A)}$ is the KV-algebra of $(\tilde{M},
\tilde{D})$ and $\tilde{W} = C^\infty (\tilde{M}, \R).$

\begin{pro} If the ``affine" manifold $(M, D)$ is complete then the cohomology class $\tilde{v}
\in H_m(\tilde{\A}, \tilde{W})$ vanishes.
\end{pro}

The proof of the above Proposition follows from the fact that, under the completeness assumption
the development map is an affine diffeomorphism from $\tilde{M}$ onto the euclidian space.

\begin{cor} Let $(M, D)$ be a compact locally flat manifold where the fundamental group
$\Gamma$ is nilpotent. Then every $\Gamma$-invariant volume form in $\tilde{M}$ is
$\tilde{\delta}$-exact, where $\tilde{\delta}$ is the boundary operator of $C(\tilde{\A},
\tilde{W})$
\end{cor}

Hint. The existence of a $\Gamma$-invariant volume form implies that $\ell(\gamma)$ is
unimodular. Following a theorem by [FGH], $(M, D)$ is complete under the assumption that $\gamma$
is nilpotent

Regarding the conjecture of Markus [CY], [FGH], we see that if $(M, D)$ is a compact
unimodular locally flat manifold and if its universal covering $(\tilde{M},
\tilde{D})$ carries a $\tilde{D}$-parallel volume form $\tilde{v}$ whose the
cohomology class $[\tilde{v}] \in H_m(\tilde{\A}, \tilde{W})$ is non-zero, then $(M,
D)$ is not complete. For instance take $M = \R^n \times \R^n - (0,0) / \Gamma \ {\rm
with} \ \Gamma = \{2^mI_{\R^n} \times 2^{-m}I_{R^n}, m \in \Z\}$, two copies of Hopf
manifolds.

\part{Groups of diffeomorphisms}

This Part III is devoted to the complexes arising from locally flat manifolds. We shall sketch
the study of left invariant affine structures on groups of diffeomorphisms and their relations
to the affine geometry of smooth manifolds.

\section{\underline{Lie groups of diffeomorphisms}}

\subsection{\underline{Generalities}}

Let $\M$ be a smooth manifold. We denote by $Diff(\M) ({\rm resp} \ Diff_o(\M))$ the group of
diffeomorphisms of $\M$ (resp the groups of diffeomorphisms which are isotopic to the
identity). In some situations the $Diff_o(\M)$ and many subgroups carry a structure of infinite
dimensional Lie group $[BA], [LJ], [RT].$

The Lie algebra of $Diff_o(\M)$ is the algebra of compactly supported smooth vector fields,
which we will denote by $\A_o(M)$.

\begin{df} A KV-structure in the $\R$-vector space $\A_o(M)$ is called a KV-structure
on $Diff_0(M)$.
\end{df}

Since we are interested in the relationship between KV-structures on $Diff_0(M)$ and
affine structures on M, there is no loss of generality in dealing with KV-structure
on $\A(M) = \X(M)$.

We already saw that there may be KV-structures on $\A(M)$ which are not defined by any locally
flat linear connection on M.

1) For instance take $M =\R$, let $\A = f \frac{d}{dx}$ and $b = g \frac {d}{dx}$, then

$$
ab = f \frac{dg}{dx} \ \frac{d}{dx}
$$

is defined by the euclidian connection on $M = \R$.

2) Take $M = \R, a = f \frac{d}{dx}, \ b = g \frac{d}{dx}$ and set

$$
ab = fg \frac{d}{dx} .
$$

We just extended a KV-structure on $X(\R)$ which is not given by a linear connection. So
it becomes interesting to ask:

(Q) what do KV-structures on $\A(M) = X(M)$ given by linear connections look like ?

To answer the question Q we shall begin by pointing out some useful tools.

\subsection{\underline{Deformations of $\Z_2$-graded KV-algebras}}

Let $\A$ be a KV-algebra. For every KV-module W we consider the semi-direct product $\A
\times W$, which is a KV-algebra with the following multiplication map :

$$
(a + w) . (a' + w') = aa' + aw' + wa' .
$$

We equip $\G = |\A| \oplus W$ with a $\Z_2$-graduation by setting

$$
\G^o = |\A| \ , \ \G^1 = W
$$

So $\G$ is a $\Z_2$-graded KV-algebra. We are interested in the KV-complex $C(\G,\G)$. The
vector space $C_q(\G,\G)$ of q-chains is $\Z_2$-graded as follows :

$$
C_q(\G,\G) = c_{q,o}(\G,\G) \oplus C_{q,1}(\G,\G)
$$

where $C_{q,j}(\G,\G)$ is the subspace of $\G^j$-valued q-multilinear maps defined in $\G$. On
the other hand, every $C_{q,j}(\G,\G^j)$ is bi-graded by those $C_{q,j}^{r,s}(\G,\G^j) =
\ho((\stackrel{r}{\otimes} |\A|) \ \otimes \ (\stackrel{s}{\otimes} W), \G^j)$ with $r + s =
q$.

Henceforth we shall set

$$
C_{q,p}^{r,s}(\G) = \ho((\stackrel{r}{\otimes} |\A|) \ \otimes \ (\stackrel{s}{\otimes} W) ,
\G^p)
$$

where $r$ and $s$ are non negative integers such that $r + s = q$ and $p \in \Z_2$.

We have $C^{r,s}_q(\G) = C^{r,s}_{q,o}(\G) \oplus C^{r,s}_{q,1}(\G).$

Suppose V is a KV-module for $\A$. Then V becomes a KV-module for $\G$ and the boundary
operator

$$
\delta: C_q(\G,V) \longrightarrow C_{q+1}(\G,V)
$$

sends $C^{r,s}(\G,V)$ into $C^{r+1,s}(\G,V)$.

Let us consider the case $V = W$. One may observe that in this case $\delta$ sends $C_{q,1}$ into
$C_{q+1,1}$, so the vector space

$$
C_{*,1}(\G,\G) = \sum_q C_{q,1}(\G,W)
$$

is a subcomplex of the chain complex $C(\G,W)$. The same claim doesn't hold for $C_{*,o}(\G,\G)
= \sum_q C_{q,o}(\G,\A)$.

\underline{Remark}. Every time there is no risk of confusion we shall simply write $a$ instead of
$(a;0)$, (resp. $w$ instead of $(0, w)$). Hence we have :

$$
\begin{array}{l}
a \in \A \subset \A \oplus W for (a,o) \in \A \times W,\\
w \in W \subset \A \oplus W for (o,w) \in \A \times W.
\end{array}
$$

Thus

$$
(a + w) (a' + w') = (a,w).(a',w') = (aa',aw' + wa').
$$

Now let us see what do elements in $J(\G)$ look like.

If $(a_o,w_o) \in J(\G)$ then for all (a,w) and (a',w') we have

$$
((a,w),(a',w'),(a_o,w_o)) = (0, 0).
$$

One sees that $a_o \in J(\A), \quad w_o \in J(W)$ and $(a,w,a_o)$ must vanish for all $(a,w) \in
\G$, so that in general we have

$$
J(\G) \subset J(\A) \oplus J(W), J(W) \subset J(\G).
$$

Therefore we deduce that the associative subalgebra $J(\G) \subset \G$ is homogeneous, viz it is
$\Z_2$-graded.

Now we are going to focus on the subcomplex

$$
\longrightarrow C^{r,s}_{q,1} \stackrel{\delta}{\longrightarrow} C^{r+1,s}_{q+1,1} (\G)
\longrightarrow.
$$

Let $\theta \in C^{1,1}_{2,1}(\G) = \ho(\A \otimes W,W).$ Then we have $\delta \theta \in
C^{2,1}_{3,1} (\G)$. Therefore, if $a, b \in \A, w \in W$ and suppose that $\theta$ is
$\delta$-closed, then we get

$$
\left \{
\begin{array}{ll}
(i) & - a \theta(b,w) + \theta(ab,w) + \theta(b,aw) = -b \theta(a,w) + \theta(ba,w) +
\theta(a,bw)\\
(ii) & - a \theta(w,b) + \theta(aw,b) + \theta(w,ab) - (\theta(w,a))b\\
& = \theta(wa,b) + \theta(a,wb) - (\theta(a,w))b.
\end{array} \right.
$$

The both conditions (i) and (ii) come from formula (25).

Let us consider $\theta \in C^{0,2}_{2,1}(\G)$.Then $\delta \theta \in C^{1,2}_{3,1}(\G)$. If
$\delta \theta = 0$, then for $a \in \A$ and $w, w' \in W$ we get

$$
\left \{
\begin{array}{ll}
(i) & -a \theta(w,w') + \theta(aw,w') + \theta(w,aw') - \theta(wa,w') = 0,\\
(ii) & \theta(w',wa) - (\theta(w',w)a = \theta(w,w'a) - (\theta(w,w'))a.
\end{array} \right.
$$

Since $\G^o$ is not a sub-KV-module for $\G$, we saw that the set of $\G^o$-valued chains is not
a sub-complex of $C(\G,\G)$. Nevertheless we are interested in some particular $\G^o$-valued
chains. Let $\theta \in C^{1,1}_{2,0}(\G)$ such that $\theta(a,w) = \theta(w,a)$ for all a $\in
\A$ and for all w $\in$ W.

If $\delta \theta$ vanishes, because $\delta \theta  \in C^{2,1}_3 \oplus C^{1,2}_3$, we have
to examine

$\begin{array}{l}
0 = \delta \theta(a,b,w) \ {\rm and} \ 0 = \delta(a,w,b) \ {\rm for} \ a , b \in
\A \ {\rm and} \ w \in W.\\
(i) \ 0 = \delta \theta(a,b,w) \Rightarrow -a \theta(b,w) + \theta(ab,w)
+ \theta(b,aw) = -b \theta(a,w) + \theta(ba,w) + \theta(a,bw)
\end{array}$

$\begin{array}{ll}
(ii) \ 0 = \delta \theta(a,w,b) \Rightarrow & -a \theta(w,b) + \theta(aw,b) + \theta(w,ab)\\
& - \theta(wa,b) - \theta(a,wb) = 0\\
(iii) \ 0 = \delta \theta(a,w,w') \Rightarrow & 0 = - \theta(w,a)w' - w \theta(a,w') +
\theta(a,w)w'
\end{array}$

By vertu of the symmetry property of $\theta(w,a)$ we see that the closeness of $\theta$ will
imply that $w \theta(a,w') = 0$ for all $a \in \A$, and $w, w' in W$. The last fact motivates
our interest in well understanding the case of $\Z_2$-graded KV-algebras $\G = \G^o \oplus
\G^1$ in which the $\G^{1,s}$ are left $\G^o$-KV-modules. In other words we have w.a = 0 for a
$\in \A$ and for all w $\in$ W.

Henceforth we assume the $\Z_2$-graded KV-algebra $\G = \A \oplus W = \G^o \oplus
\G^1$ to satisfy the condition $W.\A = (0)$. Let $\theta \in C^{0,2}_{2,1}$ be a
cocycle, then $0 = \delta \theta(a,w,w')$ will imply the relation

$$
a \theta(w,w') = \theta(aw,w') + \theta(w,aw'). \leqno (29)
$$

Thus if we equip $W = \G^1$ with the multiplication defined as follows

$$
w.w' = \theta(w,w'),
$$

then the left action of $\G^o$ look like derivation action in $(\G^1,\theta)$, viz elements in
$\G^o$ are infinitesimal automorphism of the algebra $(\G^1, \theta)$.

Now let $\theta \in C^{1,1}_{2,0}(\G)$ subject to satisfy the symmetric property

$$
\theta(w,a) = \theta(a,w).
$$

If $\delta \theta = 0$ then for all a and b in $\G^o$ and for all w in $\G'$ one has

$$
0 = -a \theta(w,b) + \theta(aw,b) + \theta(w,ab) = 0 \leqno (30)
$$

As regards to the question (Q) we shall see that the relations (30) and (31) provide useful
tools. We are interested in some special fixed $\Z_2$-graded KV-algebra $\G = \G^o \oplus \G^1$
in which $\G^1$ is a left $\G^o$-KV-module.

So the multiplication map in $\G$ is given as

$$
(a,w)(a',w') = (aa',aw'). \leqno (31)
$$

\underline{Special deformations of (31)}. Now we shall focus on the multiplication maps in $\G$
which have the following form

$$
(a,w)(a',w') = (aa', aw' + \theta(w,w')) \leqno (32)
$$

with $\theta \in C^{0,2}_{2,1}(\G)$.

\begin{df} Let $\G^o \oplus \G^1$ be a $\Z_2$-graded KV-algebra such that $\G^1$ is a left
$\G^o$-KV-module; a chain $\theta \in C^{0,2}_{2,1}(\G)$ is called a KV-chain if and only if
$(w,w',w'')_{\theta} - (w',w,w'')_{\theta} = 0$ for all w, w' and w'' in $\G^1$, where
$(w,w',w'')_{\theta}$ denotes $\theta(w,\theta(w',w'')) - \theta(\theta(w,w'), w'')$.
\end{df}

\begin{pro} The following claims are equivalent.

1) The multiplication map (32) defines a KV-algebra structure in $\G$.

2) $\theta \in C^{0,2}_{2,1}(\G)$ is a KV-cocycle.
\end{pro}

The proof consists in a direct verification using (29).

Of course a KV-chain $\theta \in C^{0,2}_{2,1}(\G)$ is called associative chain (resp.
commutative chain) if and only if the algebra $(\G^1, \theta)$ is associative (resp.
commutative). It is easy to verify that a commutative KV-chain is a associative chain.

It is to be noticed that if $\theta \in C^{0,2}_{2,1}(\G)$ is a KV-cocycle of the
$\Z_2$-graded KV-algebra $\G^o \oplus \G^1$ then the KV-algebra $\G_{\theta}$, given
by (32) is non $\Z_2$-graded KV-algebra.

Using the relation (29) and (30) we are in position to answer partially the question (Q).

Indeed, given a fixed $\Z_2$-graded KV-algebra in which $\G^1$ is left
$\G^o$-KV-module, we saw that the cocycles $\theta \in C^{0,2}_{2,1}(\G)$ will
satisfy the identity (29). At the another side, given a cocycle $\psi \in
C^{1,1}_{2,0}(\G)$, then $\psi$ will satisfy the identity (30) if and only if $\psi$
is symmetric, viz $\psi (a,w) = \psi (w,a)$ for all $a \in \G^o$ and for all $w \in
\G^1$.

Henceforth the action of $\G^o$ in $\G^1$ is faithful. We have the following
definition

\begin{df} A pair of cocycle $(\theta, \psi) \in C^{0,2}_{2,1}(\G) \times C^{1,1}_{2,0}(\G)$ is
called connection like pair if the following properties hold:

$(c_1) \ \psi \in C^{1,1}_{2,0}(\G)$ is symmetric, viz $\psi(a,w) = \psi(w,a), a \in \G^o, w \in
\G^1,$

$(c_2) \ \theta$ is a KV-cocycle,

$(c_3) \ \theta$ and $\psi$ are related as $\psi(\theta(w,w'),a) = \psi(w,\psi(w',a)); a \in
\G^o, w \ {\rm and} \  w' \ in \ \G^1$.
\end{df}

\underline{Example}. Let (M,D) be locally flat smooth manifold and let $\A = (X(M),D)$ be its
canonical KV-algebra; the vector space $W = C^{\infty} (M, \R)$ is a left KV-module for $\A$
under the Lie derivation, so that

$$
\G = \A \oplus W
$$

is a $\Z_2$-graded KV-algebra. The multiplication map is

$$
\begin{array}{ll}
(a,w)(a',w') & = (D_aa',L_aw')\\
& = (aa', D_aw')
\end{array}
$$

Let $\theta \in C^{0,2}_{2,1}(\G)$ be the usual multiplication of two smooth functions

$$
\theta(w,w') = ww'
$$

and let $\psi \in C^{1,1}_{2,0}(\G)$ be multiplication of vector fields by smooth functions

$$
\psi(w,a) = \psi(a,w) = w.a
$$

Both $\theta$ and $\psi$ are cocycles of the chain complex $C(\G,\G)$; they are
related following $(c_2)$, viz $(ww')a = w(w'a)$. On the other hand the action
$(a,w) \longrightarrow aw = \ \langle \ dw,a \ \rangle$ \ of $\A$ in W is faithful;
the cocycle $\theta$ is a KV-cocycle.

We call a connectionlike pair $(\theta, \psi)$ associative (resp. commutative)
whenever the KV-cocycle $\theta$ is associative (resp. commutative). In the example
we just set $\theta$ of order $\leq 0$ and $\psi$ is of order $\ \leq \ 1$.

Now, given a connectionlike pair $(\theta, \psi)$, we deduce the deformation $\G_{\theta}$ of
$\G$ which has the following multiplication map:

$$
(a,w)(a',w') = (aa',aw' + \theta (w,w')).
$$

So that the pair $(\G_{\theta}, \psi)$ will be called a connectionlike deformation of the
$\Z_2$-graded KV-algebra $\G = \G^o \oplus \G^1$.

\underline{Formally}: given a connectionlike deformation $(\G_{\theta}, \psi)$ of a
$\Z_2$-graded algebra $\G = \G^o + \G^1$ ; \underline{then}:

$(*)_1: \G^o$ is taken to be the space of ``vector fields" on some ``manifold".

$(*)_2: \G^1$ is taken to the space of ``smooth functions" on that manifold.

$(*)_3$: the multiplication in $\G^o$ is taken to be a linear connection.

$(*)_4$: the left action of $\G^o$ in $\G^1$ is the Lie derivation: (the analogue the bracket of
vector fields is $[a,a'] = aa' - a'a)$.

$(*)_5$: the cocycle $\psi$ is the multiplication of vector fields by smooth
functions and the identities (29) and (30) are well understood.

\begin{Th} Let $(\G_{\theta}, \psi)$ be a connectionlike deformation of the $\Z_2$-graded
KV-algebra $\G^o \oplus \G^1$ such that $\theta \in C^{0,2}_{2,1}(\G) \setminus \{0\}$. Then the
one parameter family $t \longrightarrow \G_t = (\G_{t \theta},\psi)$ is an one parameter family
of non trivial connectionlike deformations of $\G^o \oplus \G^1$.
\end{Th}

\underline{Hint}. For every $t \in \R, t{\theta}$ is a KV-cocycle and the pair $(t
\theta , \psi)$ are related as in Definition 10.3. Of course $\theta \in
C^{0,2}_{2,1}(\G) - \{0\}$ is non $\delta$-exact; so that the deformation $\G_t$ is
non trivial.

To end considerations regarding the question (Q) it is to be noticed that
connectionlike deformations of $\Z_2$-graded KV-algebras are related to special
2-cocycles. Indeed 2-chains $c \in C_2(\G,\G)$ have six components, said
$c^{r,s}_{2,j}$ with $r \ {\rm and}\ s \ {\rm in} \ \{0,1,2\} \ {\rm and}\ j \ {\rm
in}\ \Z_2$.

\begin{df} 1) KV-2-chains of $\G$ are those elements $c \in C_2(\G,\G)$ the components $\theta =
c^{0,2}_{2,1}$ of which are KV-multiplications in $\G^1$, viz

$$
(w,w',w'')_{\theta} = (w',w,w'')_{\theta}
$$

for all w, w' and w'' in $\G^1$.

2) Connectionlike KV-2-chains of $\G$ are those KV-chains $c \in C_2(\G,\G)$ with only two
components $\theta = c^{0,2}_{2,1} \ {\rm and} \ \psi = c^{1,1}_{2,0}$ such that
$c^{1,1}_{2,0}$ is regular and symmetric, viz $\psi(a,w) = \psi(w,a)$ for all $a \in \G^o$ and
for all $w \in \G^1$.
\end{df}

\underline{Remark}. Let $c \in C_1(\G,\G)$ and let us consider the $\G^o$-component
of c, denoted by $\xi : \G \longmapsto \G^o. \ {\rm Set} \ \psi = \delta \xi$, and
apply $\psi$ to $((a,o),(o,w)) = (a,w)$. It will be symmetric w.r.t (a,o) and (o,w)
iff $\psi(a,w) = 0$ so that if a connectionlike cocycle $(\theta,\psi)$ is exact
then it vanishes identically.

\begin{Th}. There is an one-to-one correspondence between the set of connectionlike KV-2-cocycles
in $C_2(\G,\G)$ and the set of connectionlike deformations of the $\Z_2$-graded KV-algebra $\G =
\G^o \oplus \G^1$.
\end{Th}

\underline{Sketch of the proof}. Let $c \in C_2(\G, \G)$ be a connectionlike cocycle
of the $\Z_2$-graded KV-algebra $\G = \G^o \oplus \G^1$. Since c have two components
only $c^{1,1}_{2,0} \ {\rm and} \ c^{0,2}_{2,1}$ let us set $\theta(w,w') =
c^{0,2}_{2,1}(w,w') \ {\rm and} \ \psi(a,w) = \psi(w,a) = c^{1,1}_{2,0}(a,w) =
c^{1,1}_{2,0}(w,a)$. Of course we write $(a,w) \simeq (a,o) + (o,w)$ and so on.

One has $\delta \theta \in C^{1,2}_{3,1}(\G,\G)$ and $\delta \psi \in C^{2,1}_{3,0}(\G,\G) +
C^{1,2}_{3,0}(\G)$.

For $(a,w)$ and $(a',w')$ in $\G = \G^o + \G^1$, we have

$$
c \ ((a,w),(a',w')) = (\psi(a,w') + \psi(w,a'), \theta(w,w')) \in \G^o \oplus \G^1.
$$

Since $\G^1$ is a left $\G^o$-KV-module, the closeness condition $\delta c = 0$ implies the
following system:

$$
\left \{
\begin{array}{l}
1) \ a \theta(w',w'') - \theta(aw',w'') - \theta(w',aw'') = 0,\\
2) \ a \psi(a',w'') - \psi(aa',w'') - \psi(a',aw'') = 0,\\
3) \ \psi(a, \theta(w',w'')) - \psi(\psi(a,w'),w'') = 0.
\end{array} \right.
$$

The identities 1) and 2) tell us that $\theta$ and $\psi$ are 2-cocycles with
properties (29) and (30) respectively. The condition 3) tells us that $\theta$ and
$\psi$ are related as in the definition 10.3.

Since c is assumed to be a KV-2-cocycle one has

$$
(w,w',w'')_{\theta} = (w',w,w'')_{\theta}
$$

for all w,w' and w'' in $\G^1$. Now we define the KV-algebra $\G_c = \G_{\theta}$ whose
multiplication map is given as

$$
(a,w).(a',w') = (aa',aw' + \theta(w,w')).
$$

The condition above implies that $(\G_c, \psi)$ is a connectionlike KV-algebra.

That ends the proof of the theorem 10.2.

\section{Left invariant KV-structures on $Diff_0(M)$}

Let G be a finite dimensional Lie group and $\G$ its Lie algebra of left invariant
vector fields. Let $KV(|\G|)$ be the set of KV-algebra structures on $|\G|$. For
every $\mu \in KV(|\G|)$, the commutator $ [a,a'] \mu = \mu(a,a') - \mu(a',a)$
defines a Lie algebra structure in $|\G|$, denoted by $\G_{\mu}$. Not every $\mu \in
KV(|\G|)$ will give rise to a left invariant linear connection on the Lie group $G$.

In order that $\mu \in KV(|\G|)$ give rise to a left invariant linear connection on
$G$, it is necessary and sufficient that $ [a,a']_{\mu} = [a,a'],$ where $[a,a']$ is
the Lie bracket of vector fields on the smooth manifold G. We denote by $KV(\G)$ the
subset of $\mu \in KV(|\G|)$ such that $[a,a']_{\mu} = [a,a']$ for all $a$ and $a'$
in $\G$. So $KV(\G)$ consists of connectionlike KV-algebra structures on $|\G|$.
Mutatis mutandis we will attempt to write the (formal) analogue to $KV(|\G|)$ and
$KV(\G)$ for Lie groups of diffeomorphisms. Let M be a smooth manifold and let
$Diff_0(M)$ be the group of those diffeomorphisms of M which are isotopic to the
identity map. It is well known that $Diff_0(M)$ carries an infinite dimensional Lie
group structure whose Lie algebra is $\A_c(M)$, the Lie algebra of compactly
supported smooth vector fields, [LJ], [RT].

Now we have seen that $\A_c(M)$ may carry a KV-algebra structure which is not given
by a locally flat linear connection on M. (See 10.1, Example 2). Nevertheless by
imitating the situation of finite dimensional Lie group we set the following
definition.

\begin{df} A left invariant KV-structure in $Diff_0(M)$ is a connectionlike KV-algebra structure
in the $\Z_2$-graded vector space $\A_c(M) \oplus C^{\infty}(M,\R)$ whose commutator
Lie algebra induces on $\A_c(M)$ the Lie algebra $(\A_c(M), [\ ,\ ])$ where $[\ ,\
]$ is the Lie bracket.
\end{df}

\underline{Example}: let $(M,D)$ be a locally flat smooth manifold and $\A =
(X(M),D)$ its KV-algebra. The subspace $\A_c(M)$ of compactly supported smooth
vector fields is a two-sided ideal in $\A$. So the $\Z_2$-graded subspace $\A_c(M)
\oplus C^{\infty}(M,\R)$ is a bilateral $\Z_2$-graded ideal in $\A \oplus
C^{\infty}(M,\R)$, the multiplication map being $(a,w).(a',w') = (aa',aw')$, with
$a,a'$ in $\A$, $w$ and $w'$ in $C^{\infty}(M,\R)$ and $aa' = D_aa', aw = D_aw$. The
datum $(M,D)$ gives rise to a linear connectionlike deformation of $\A \oplus
C^{\infty}(M,\R)$ :
$$
(a,w)(a',w') = (aa',aw' + ww'),
$$
where $ww'$ is the usual product of two functions. Of course the connectionlike
KV-algebra above induces $\A_c = (X_c(M),D)$. Thus $Diff_0(M)$ carries a left
invariant locally flat structure (or left invariant KV-structure).

If $M$ is without boundary, then $Diff_0(M)$ is a simple Lie group, subject to carry
left invariant KV-structures. This is in contrast with the case of finite
dimensional semi simple Lie groups $[NB3]$. So, starting with a locally flat
manifold $(M,D)$, we know that $Diff_0(M)$ will carry a left invariant KV-structure.
An interesting problem consists of the search of finite dimensional Lie subgroups $G
\subset Diff_0(M)$ on which the left invariant KV-structures of $Diff_0(M)$ induce
KV-structures. In particular $Diff_D(M)$, the group of $D$-preserving
diffeomorphisms deserve some attention. Inversely, it would be interesting to
characterize those left invariant KV-structure on $Diff_0(M)$ which really induce
locally flat structures in $M$. For some manifolds we know that no much structure
does exist. For instance if $M$ is the three sphere $S^3$.

\bibliographystyle{plain}

\end{document}